\documentclass[review]{elsarticle}

\usepackage{lineno,hyperref}
\journal{J. Math. Anal. Appl.}
\usepackage{amscd}
\usepackage{mathrsfs}
\usepackage{amsfonts}
\usepackage{graphicx}
\usepackage{amsmath,amscd,stmaryrd,amssymb,array, amsfonts,amsmath,amssymb}
\usepackage{enumitem}

\usepackage{caption}
\usepackage{subfigure}

\numberwithin{figure}{section}
 \numberwithin{equation}{section}

\newtheorem{theorem}{Theorem}[section]
\newtheorem{proposition}[theorem]{Proposition}
\newtheorem{definition}[theorem]{Definition}

\newtheorem{lemma}[theorem]{Lemma}
\newtheorem{remark}[theorem]{Remark}

\newcommand{\bi}{{\mathbf{i}}}

\newcommand{\bbA}{{\mathbb A}}

\newcommand{\bbC}{{\mathbb C}}

\newcommand{\bbP}{{\mathbb P}}

\newcommand{\bbS}{{\mathbb S}}

\newcommand{\cD}{{\mathcal D}}

\newcommand{\sL}{{\mathscr L}}

\def\be{\begin{equation}}
\def\ee{\end{equation}}
\def\bes{\begin{equation*}}
\def\ees{\end{equation*}}
\def\bsp{\begin{split}}
\def\esp{\end{split}}

\def\ba{\begin{array}}
\def\ea{\end{array}}
\def\benu{\begin{enumerate}}
\def\eenu{\end{enumerate}}
\def\bt{\begin{theorem}}
\def\et{\end{theorem}}
\def\bp{\begin{proposition}}
\def\ep{\end{proposition}}
\def\bl{\begin{lemma}}
\def\el{\end{lemma}}
\def\br{\begin{remark}}
\def\er{\end{remark}}
\def\bd{\begin{definition}}
\def\ed{\end{definition}}

\def\b{\beta}

\def\de{\delta}

\def\pa{\partial}

\def\lam{\lambda}

\def\o{\stackrel{\circ}}
\def\ve{\varepsilon}
\def\sig{\sigma}

\def\gam{\gamma}

\def\a{\alpha}

\def\.{\cdot}

\def\bbA{\mathbb{A}}
\def\bbC{\mathbb{C}}

\def\R{\mathbb{R}}

\def\bbX{\mathbb{X}}

\def\A{\forall}
\def\ol{\overline}

\def\Cap{\bigcap}

\def\ra{\rightarrow}

\def\~{\tilde}
\def\8{\infty}
\def\X{\times}
\def\({\left(}
\def\){\right)}

\def\mb{\mbox}
\def\emp{\emptyset}
\def\sm{\setminus}

\def\bx{\blacksquare}

\def\Hs{\hspace{1cm}}\def\hs{\hspace{0.5cm}}
\def\Vs{\vskip8pt}\def\vs{\vskip4pt}

\def\({\left(}\def\){\right)}

\begin{document}

\begin{frontmatter}

\title{A Dynamical Approach to the Perron-Frobenius Theory and Generalized Krein-Rutman\\ Type Theorems$^\dag$\footnote{$\dag$This work was supported by the National Natural Science Foundation of China [11871368]}}
%\tnotetext[]{{\em 2010 Mathematics Subject Classification}. Primary  47B65, 47A75, 15A18, 34A30.}
%% Group authors per affiliation:
%\author{Luyan Zhou\fnref{myfootnote}}
%\address{Radarweg 29, Amsterdam}
%\fntext[myfootnote]{Since 1880.}

%% or include affiliations in footnotes:

\author[mymainaddress]{Desheng Li\,\corref{mycorrespondingauthor}}
\cortext[mycorrespondingauthor]{Corresponding author}
\ead{lidsmath@tju.edu.cn}

\author[mymainaddress]{Mo Jia\,}
\ead{jiamomath@tju.edu.cn, jiamomath@sina.com}
%\address[mymainaddress]{Center for Applied Mathematics,  Tianjin University, Tianjin 300072,  China}
\address[mymainaddress]{School of Mathematics,  Tianjin University, Tianjin 300072,  China}

\begin{abstract}
We present  a dynamical approach to the classical Perron-Frobenius theory by using  some elementary  knowledge on  linear ODEs. It is completely self-contained and significantly different from those  in the literature.
As a result, we develop a complex version of the Perron-Frobenius  theory and prove some   generalized Krein-Rutman type theorems.
\end{abstract}

\begin{keyword}
%\texttt{elsarticle.cls}\sep \LaTeX\sep Elsevier \sep template
%sectorial operator,
Perron-Frobenius theory\sep Krein-Rutman theorem\sep dynamical approach,  rotational strong positivity, weak irreducibility.
\vs\MSC[2010] 34A30\sep 47B65\sep 47A75\sep 15A18.
\end{keyword}
\end{frontmatter}

%\linenumbers

%\newpage
\tableofcontents
%\newpage

\section{Introduction}

 In the pioneering work \cite{Perron2} Perron published  his celebrated theorem asserting  that  the spectral radius of a positive matrix is an eigenvalue of the matrix which corresponds to a positive eigenvector. This result was later
 extended by  Frobenius to  nonnegative irreducible matrices \cite{Frobenius}, which is now known as the Perron-Frobenius ({\em PF}\, in short) theorem. Since in practice it is often the case that a matrix has  nonnegative entries, the theorem has found an astonishing wide range of applications in diverse areas; % such as the numerical analysis, wireless networks,  economy, epidemiology, probability, dynamical systems theory, lower-dimensional topology, etc.;
 see e.g. MacCluer \cite{MacCluer} for a review on this topic. This greatly  stimulated  mathematicians to seek for  various  proofs from different perspectives,  aiming at  opening  up new ways to this fruitful area.   Many elegant  proofs of the theorem were summarized in \cite{MacCluer}, including the Perron's original proofs  \cite{Perron1,Perron2} based on a technical use of the Cramer's rule or the resolvent of matrices, Wielandt's 1950 proof \cite{Wiel}  built on Frobenius's mini-max  idea, Birkhoff's proof  via  the Jordan canonical form \cite{Birkhoff}, Karlin's proof using the complex variable theorem about power series with positive coefficients \cite{Karlin}. Other   proofs  based on geometric methods and  fixed-point theorems can be found in \cite{B2,Koh,Pull,Same}, etc.

  As we have seen in the literature, the new ideas and techniques involved in different proofs have led to many  important   extensions of the  PF theorem. An outstanding extension  is the famous Krein-Rutman  ({\em KR} in short)  theorem established in 1948 \cite{Krein}  dealing with positive compact linear operators, which was obtained by a substantial  use of the  Schauder fixed-point theorem. % and the spectral theory about  compact operators.
  This infinite-dimensional extension tremendously expanded the scope of applications of the PF  theory, and therefore aroused a great interest   in further studies from different perspectives;
   see e.g.  \cite{Arendt,Kry,Kras,Nuss,Nussb,Schab}.

 Another remarkable  one is the  nonlinear PF theory,  whose development was originally   due to a   simple-looking %geometric and topological
 observation   made by  Birkhoff \cite{B2} and   Samelson \cite{Same} (who remarked that one can
use Hilbert's projective metric and the contraction mapping principle to
prove PF type theorems for nonnegative matrices). The interested reader may consult % Lemmens and Nussbaum
\cite{LN} for a finite-dimensional nonlinear PF theory, and \cite{Maha,Paret,Ogi} etc. for infinite dimensional results in this line.
In recent years there even  appeared some important extensions of the PF theorem to   positive tensors   \cite{CPZ,YY}.

%     A dynamical systems proof  was given in Alikakos and Fusco \cite{Alikakos}. The basic idea in \cite{Alikakos} is to  consider the  discrete dynamical system of the induced operator on the projective space, and the existence of a positive eigenvector of the original operator reduces to showing that the $\omega$-limit set of the induced operator consists of a single equilibrium. %This is essentially different from that of ours in this present work.
  % In  \cite{CQ} Chang,  Pearson and Zhang established  some Krein-Rutman type theorems for tensors. There are also elegant  works  in the literature concerning nonlinear extensions of the KR theorem; see e.g. Chang \cite{Chang1}, Mahadevan \cite{Maha} and Mallet-Paret and Nussbaum \cite{Paret}, etc.

The main purpose of this  paper is twofold. One is of a methodological sense. Specifically, we  present  a   self-contained  new dynamical  approach towards the PF   theory by using only some elementary knowledge on linear ODEs; this approach is significantly different from those in the literature and  gives us a  dynamical picture  of the theory.
The other is to develop a complex version of the theory and prove some  generalized KR type theorems for  bounded real  operators.

This work will be carried out  as follows. In Section \ref{s:2} we introduce some basic notions and notations concerning cones and operators.  In particular, %some basic   notions used in the real PF theory seem to be somewhat restrictive. As one of the main ingredients of  this section,
we introduce the  notion of  {\em rotational strong positivity} for complex  operators to replace the one of {strong positivity} which  plays a crucial role  in the real  KR theorem. Some simple  facts used throughout the paper  will also be presented.
In Section \ref{s:3} we prove two technical lemmas. One concerns a fundamental projective  property  of cones, and the other relates to  the asymptotic behaviour  of linear ODEs in finite-dimensional complex Banach spaces.

Sections \ref{s:4} and \ref{s:5} are concerned with a finite-dimensional complex version of the PF theory. Specifically, we develop in a self-contained manner an elementary dynamical approach towards the theory  and prove in Section \ref{s:4}   a generalized   finite-dimensional KR type theorem  for complex  operators  under weaker assumptions (Theorem \ref{t:4.1}). This result   immediately yields  a refined  complex version of the  KR theorem  as long as we  impose on the operator the rotational strong positivity assumption, which forms the main ingredients of    Section \ref{s:5}.

In Section \ref{s:6} we extend the finite-dimensional  results in previous sections to infinite-dimensional spaces
 by a simple use of direct sum decompositions of the spaces which reduces infinite-dimensional problems to finite-dimensional ones.

 Section \ref{s:7} is devoted to real bounded linear  operators. By considering the complexifications of real   operators and applying the results for complex operators, we obtain a generalized  KR type  theorem (Theorem \ref{t:7.4}) and give simple proofs for some stronger  versions of the classical KR theorem. We mention that a nice  dynamical systems proof for a major part of the KR theorem can be found  in Alikakos and Fusco \cite{Alikakos}. But our  proof here is of a  completely   different nature from the one given  in \cite{Alikakos}.

\section{Preliminaries}\label{s:2}
%In this  section  we first introduce some notions and notations that will be used throughout the paper. Then we give  some basic facts concerning cones and operators in Banach spaces. In particular, we prove a  projective  property of cones with respect to direct sum  decompositions of spaces.

\vs
Throughout the  paper $\R$ and    $\bbC$ denote respectively  the fields of real and complex numbers, and  $\mathbf{i}$  the unit imaginary number in $\bbC$. Let $\bbS_r=\{w\in\bbC:\,\,|w|=r\}$ ($r>0$), and  $\R_+=[0,\8)$.

Let   $E$ be  a Banach space (real or complex) with norm $\|\.\|$.
Given a set   $M\subset E$, the {\em interior} and {\em closure } of $M$  are denoted respectively by  $\o M$ and $\ol M$. When we need to emphasize in which space the interior and closure are taken, we may use the notations $\mb{int}_E M$ and $\mb{Cl}_E M$ in place of  $\o M$ and $\ol M$, respectively.
Given   $x\in E$, set $d(x,M)=\inf_{y\in M}\|x-y\|.$

Denote by   $\sL(E)$  the space of bounded linear  operators on $E$.

\subsection{Basic knowledge in the spectral theory  of operators}\label{s:2.3}
\vs Let $E$ be  a complex Banach space, and  $A\in\sL(E)$.  Denote by $\sig(A)$ and $\rho(A)$ the {\em spectrum} and {\em resolvent set} of $A$, respectively.  For $\lam\in\rho(A)$, let $R_\lam(A):=(\lam-A)^{-1}$ be the resolvent of $A$.

 Given  $\mu\in\sig( A)$, set
\be\label{e:es}
\mb{GE}_\mu( A)=\{\xi\in E:\,\,( A-\mu)^{j}\xi=0\mb{ for some }j\geq 1\}.
\ee
%called the {\em generalized eigenspace} of $ A$ pertaining  to $\mu$.
It is trivial to check that $A$ leaves $\mb{GE}_\mu( A)$ invariant.
% i.e., $A(\mb{GE}_\mu( A))\subset \mb{GE}_\mu(A)$.
If $\mb{GE}_\mu(A)\ne\{0\}$, then  we call each  $\xi\in \mb{GE}_\mu( A)\sm\{0\}$  a {\em generalized eigenvector}  of $ A$. Clearly  for such a  $\xi$,  there is an  integer $k\geq 1$ such that
 %for each $\xi\in \mb{GE}_\mu( A)$ with $k=\mb{rank}(\xi)$, we have
\be\label{e:2.5}( A-\mu)^{j}\xi\ne0\,\,\,(0\leq j\leq k-1),\hs ( A-\mu)^{k}\xi=0.\ee
(Here we assign $(A-\mu)^0=I$, the identity mapping on $E$.) It is convenient to  call the number $k$  in \eqref{e:2.5}  the {\em rank} of $\xi$, denoted by  $\mb{rank}(\xi)$.
Note that an  {\em  eigenvector } $w$  is a generalized eigenvector with $\mb{rank}(w)=1$.

Denote by  $\sig_e(A)$     the {\em essential spectrum}   of $A$ in the terminology of Browder \cite[pp. 107-108, Def. 11]{Browder}. Then  each $\mu\in\sig(A)\sm\sig_e(A)$ is isolated in $\sig(A)$ with $\mb{GE}_\mu(A)$ being a finite-dimensional subspace of $E$; see \cite[pp. 108]{Browder}.

If $A\in\sL(E)$, we define   the {\em spectral radius} $r(A)$ and  {\em essential spectral radius} $r_e(A)$ as
$$r(A)=\sup\{|\mu|:\,\,\mu\in\sig(A)\},\hs r_e(A)=\sup\{|\mu|:\,\,\mu\in\sig_e(A)\}.$$
%$r_\sig$ and $r_e$ are called the {\em spectrum radius} and the {\em essential spectrum radius} of $A$, respectively.
It is basic knowledge   that $r(A)=\lim_{k\ra\8}\|A^k\|^{1/k}$.
%\be\label{e:r}r(A)=\lim_{k\ra\8}\|A^k\|^{1/k}.\ee

%The following fundamental  fact  will   be used   in Section \ref{s:8} when dealing with unbounded operators.
 % \bl\label{l:2.6}  Let $\mu\in\sig(A)$, and  $\lam\in \rho(A)$. Then
%\be\label{e:2.6}R_\lam(A)\mb{\em GE}_\mu( A)= \mb{\em GE}_\mu(A)=\mb{\em GE}_{(\lam-\mu)^{-1}}(R_\lam(A)).\ee\el
%{\bf Proof.} Let $\xi\in \mb{GE}_\mu(A)$. Then $(A-\mu)^k\xi=0$, where  $k=\mb{rank}(\xi)\geq 1$. Hence
%$$(A-\mu)^k\(R_\lam(A)\xi\)=R_\lam(A)(A-\mu)^k\xi=R_\lam(A)0=0.$$
%It follows  that $R_\lam(A)\xi\in \mb{GE}_\mu(A)$. Simple calculations also yield
%\be\label{e:2.8}
%(\lam-A)^k\(R_\lam(A)-\lam_\mu\)^k\xi=\lam_\mu^k(A-\mu)^k\xi=0,
%\ee
%where $\lam_\mu=(\lam-\mu)^{-1}$. Since $\lam\in\rho(A)$, \eqref{e:2.8}  implies   that $\(R_\lam(A)-\lam_\mu\)^k\xi=0$. Therefore $\xi\in \mb{GE}_{\lam_\mu}(R_\lam(A))$.
%In conclusion, we have $$R_\lam(A)\mb{GE}_\mu( A)\subset \mb{GE}_\mu(A)\subset\mb{GE}_{\lam_\mu}(R_\lam(A)).$$

%The verification of the inverse inclusions is similar. We omit it. $\bx$

\subsection{Cones}
\vs Let $E$ be a (real or complex) Banach space.
A {\em wedge} in $E$ is a   closed subset  $ P\subset E$ satisfying  that   $t P\subset  P$ for all $t\geq 0$.

A convex wedge $P$ with $ P\cap (- P)=\{0\}$ is called a {\em cone}.  It is  basic knowledge that if $P$ is  a cone then
 $P+P\subset P.$
 %One can also easily verify  that if $E_0$ is a closed subspace of $E$, then $P_0:=P\cap E_0$ is a cone in $E_0$, provided  that $P_0\ne\{0\}$.

Let $P$ be a cone in $E$. $ P$ is said to be   {\em total} (resp., {\em solid}),  if $\ol{P-P}=E$ (resp., $\o P\ne\emp$).
We remark that the notion of totalness of a cone $P$ is  a bit  restrictive  in the complex case, because in such a case  $P-P$ may fail to be a subspace of $E$.
Instead,  we introduce a  weaker notion  which works equal well as totalness.
\bd We say that $P$ is {full}, if there is no proper closed subspace $E_0$ of $E$ such that  $P\subset E_0$. \ed

%In other words,  $P$ is called {full} if $E_0$ is a closed subspace of $E$ and $P\subset E_0$, then $E_0=E$.
It is easy to deduce  that  a solid  cone  is total, and a total cone is full.
However, in the complex case a full cone may not be  total.

\br {\em Some authors  use the word {\em full} with
a different meaning in the context of ordered Banach spaces. For instance, in \cite{KG}  a subset $M$ of  an ordered real  Banach space  is
called {\em full} if the order interval $[a, b]$ is contained in $M$ for all $a,b\in M$.}
\er

\noindent{\bf Example 2.1.} Let
$$
P_1=\{(x,0)\in \bbC^2:\,\,x\in\R_+\},\hs P_2=\{(0,y)\in \bbC^2:\,\,y\in \R_+\}.
$$
Put  $P=P_1+\bi P_2$. Then $P$ is a cone in $\bbC^2$. Note that
$$\ba{ll}
P-P&=(P_1+\bi P_2)-(P_1+\bi P_2)=(P_1-P_1)+\bi (P_2-P_2)=X_1+\bi X_2,\ea
$$
where $X_1=\{(x,0):\,\,x\in\R\},$ and $X_2=\{(0,y):\,\,y\in\R\}.$ Clearly $\ol{P-P}=P-P\ne \bbC^2$.
Hence  $P$ is not total.

Nevertheless, $P$  is full since it contains the basis $\{(1,0),\,(0,\bi)\}$ of $\bbC^2$.

\br\label{r:2.9}{\em If  $E$ is  a real Banach space, one easily verifies that $P$   is full if and only if it is total.}
\er

The following simple fact will  be frequently used throughout the paper.
\bl\label{l:2.4} Let $ E_0$ be a finite-dimensional subspace of $E$ with $P\cap E_0=\{0\}$, and let $\{x_k\}_{k\geq 1}\subset P$ be a sequence. Suppose $d(x_k,E_0)\ra 0$ as $k\ra\8$. Then $$\mb{$x_k\ra0$\,\, as $k\ra\8$}.$$

\el
\noindent{\bf Proof.} We include a proof for completeness. Suppose on the contrary that there is a subsequence of $\{x_k\}_{k\geq 1}$, still denoted by $\{x_k\}_{k\geq 1}$, such that $\|x_k\|\geq\de>0$ for all $k$. Let $y_k=x_k/\|x_k\|$. Then $\|y_k\|=1$. We observe that
\be\label{e:2.7}\ba{ll}
d(y_k,E_0)=d\(\frac{x_k}{\|x_k\|},\,E_0\)=\frac{1}{\|x_k\|}d(x_k,E_0)\ra 0\ea
\ee
as $k\ra\8$. (The second equality in \eqref{e:2.7} is due to the fact that $E_0$ is a subspace.)
Thus one can find a bounded sequence $\{z_k\}_{k\geq 1}\subset E_0$ such that $\|y_k-z_k\|\ra 0$. Since $E_0$ is  finite-dimensional, up to a subsequence it can be assumed that  $z_k\ra z_0\in E_0$. Then  $y_k\ra z_0$ as well. Hence $\|z_0\|=1$. On the other hand, because  $y_k\in P$ for all $k$, the closedness of $P$  implies that $z_0\in P$. Therefore  $P\cap E_0\ne\{0\}$, a contradiction. $\bx$

%The main purpose in this subsection is to introduce some notions and notations concerning linear operators.
\subsection{Rotational strong positivity of complex operators}\vs

Recall that in  case  $E$ is a real Banach space and $P$ is   a solid cone in $E$, an operator $A\in\sL(E)$ is called {\em strongly positive} if $A(P\sm\{0\})\subset \o P$.
This notion  plays a crucial role in the real KR theorem. However, it is not suitable  for complex operators.  This can be seen from the following simple observation.

Let $E$ be a finite-dimensional complex Banach space, and $P$ a solid cone in $E$. Let $A\in\sL(E)$. Suppose that  $A$ is strongly positive in the terminology as above. Then as in Theorem \ref{t:5.2} we deduce  that $r=r(A)$ is an eigenvalue of $A$ with a corresponding eigenvector $w\in\o P$. On the other hand, noticing that $-w\not\in P$, one easily deduces that there is a $z_1\in\bbS_1$ such that $w_1:=z_1w\in\pa P$. Because $w_1$ is also an eigenvector corresponding to $r$, we then have $Aw_1=rw_1\in\pa P$, which  contradicts the strong positivity hypothesis.

Let us also mention   that the complexification of a strongly positive real operator may fail to be strongly positive (see Section \ref{s:7} for detail).
%To obtain  refined information concerning principal eigenvalues, here

 In what follows we introduce an  alternative  notion, called  {\em rotational strong positivity}, for complex operators.

\bd\label{d:2.1} Let $E$ be a complex Banach space, and  $P$ a solid cone. An operator $A\in\sL(E)$ is called  {rotationally strongly positive}, if\,
$\bbS_1(A x)\cap \o  P\ne\emp$ for all $x\in P\sm\{0\},$ where $\bbS_1(Ax)=\{zAx:\,\,z\in\bbS_1\}$.
\ed

 \noindent{\bf Example 2.2.} The complexification of a strongly positive real operator is rotational  strongly positive; see Section \ref{s:7}.
 \Vs

\noindent{\bf Example 2.3.} Let $|\.|$ be the usual norm of $\bbC$, and set   $|z|'=|x|+|y|$\, for $z=x+\bi y\in\bbC$.
Note that if $$z_1,z_2\in \bbC_+= \{x+\bi y\in \bbC:\,\,x,y\geq 0\},$$ then $|z_1+z_2|'=|z_1|'+|z_2|'$. It is trivial to check that
$$
|z|\leq |z|'\leq \sqrt 2|z|,\Hs\A\,z\in\bbC.
$$

Let $E:=\bbC^2=\{(z,w)^T:\,\,z,w\in\bbC\}$, where $(z,w)^T$ denotes the transpose of $(z,w)$. Define
\be\label{e:2.4}
P=\{(z,w)^T \in \bbC_+\X\bbC:\,\,|w|\leq |z|'\}.
%\(\begin{matrix} \end{matrix}\)
\ee
  We claim  that $P$ is a cone in $E$. To see this, it suffices to check that $P$ is convex and  $P\cap(-P)=\{0\}$.
Let $u_i=(z_i,w_i)^T \in P$, $i=1,2$.
Then   for every
$\theta\in[0,1]$,
$$\ba{ll}
|\theta w_1+(1-\theta)w_2|&\leq \theta |w_1|+(1-\theta)|w_2|\\[1ex]
&\leq \theta\,|z_1|'+ {(1-\theta)}|z_2|'=|\theta z_1+(1-\theta)z_2|'.
%& \leq \mb{Re}\,(\theta z_1+(1-\theta)z_2)+\mb{Im}\,(\theta z_1+(1-\theta)z_2).
\ea
$$
Hence by definition we have $\theta u_1+(1-\theta)u_2\in P$.  This verifies the convexity of $P$.
Now assume that  $u=(z,w)^T \in P\cap(-P)$. We show that $u=0$, thus proving what we desired. Suppose on the contrary that $u\ne 0$. Then either $z\ne 0$, or $w\ne 0$. If $z=0$, since $u\in P$, by the definition of $P$ one would have $w=0$. This contradicts the hypothesis $u\ne 0$.  Therefore we see that $z\ne 0$. It then follows by $u\in-P$ that $-z\in \bbC_+$. Since $z\in \bbC_+$, this  leads to a contradiction.

Note that $P$ is a solid cone in $E$ with $$
\o P=\{(z,w)^T \in P:\,\,z\in  \o\bbC_+,\,\,|w|<|z|'\},
%\(\begin{matrix} \end{matrix}\)
$$
where  $  \o\bbC_+=\{z=x+\bi y:\,\,x>0,y>0\}$.

   Let
$
A=\(\begin{matrix}a&b\\ c&d \end{matrix}\),
$
where $a>0$ is a   real number, and $b,c,d\in \bbC$. Then
$$
Au=\(az+bw,\, cz+dw\)^T ,\Hs \A\,u=\(z,w\)^T\in E.$$
In the following  we show that $A$ is rotationally  strongly positive if
 \be\label{e:5.6d}|b|<1/2,\hs 0<|c|+|d|\leq \sqrt2(a-2|b|)/4.\ee

Let $u=(z,w)^T \in P$, $u\ne 0$. Then $z\ne0$. We may write $z=re^{\bi \b}$ ($r>0$, and $\,0\leq \b\leq \pi/2$).
%Thus
%\be\label{e:5.5}|z|'=r(\cos \b+\sin \b)=r\sqrt{1+\sin2\b}.\ee
Let $\lam=e^{\bi (\pi/4-\b)}$. Clearly $\lam z=re^{\bi \pi/4}\in \o \bbC_+$. We observe  that
\be\label{e:5.6}\ba{ll}
|\lam(cz+dw)|&=|cz+dw|\leq |c||z|+|d||w|\\[1ex]
&\leq |c|| z|+|d|\,|z|'=(|c|+|d|)|z|'.
\ea
\ee
Note that $|z|'\leq \sqrt2|z|=\sqrt2|\lam z|\leq \sqrt2|\lam z|'$. Therefore
\be\label{e:5.6b}\ba{ll}
|\lam(cz+dw)|\leq \sqrt2(|c|+|d|)|\lam z|'<2{\sqrt2}(|c|+|d|)|\lam z|'.\ea
\ee

Let $\~w=az+bw$. Then
$$\ba{ll}
|\lam \~w|'&\geq |\lam a z|'-|\lam bw|'\geq |\lam a z|'-\sqrt2|\lam bw|\\[1ex]
&= |\lam a z|'-\sqrt2|b||w|\geq |\lam a z|'-\sqrt2|b||z|'\\[1ex]
&\geq a|\lam z|'-2|b||\lam z|'=(a-2|b|)|\lam z|'.
\ea
$$
Combining this with \eqref{e:5.6b} it yields
\be\label{e:5.6c}\ba{ll}
|\lam(cz+dw)|< \frac{2\sqrt2(|c|+|d|)}{(a-2|b|)}|\lam\~w|'\leq (\mb{by }\eqref{e:5.6d})\leq |\lam\~w|'.\ea
\ee

Recall that $\lam z=re^{\bi\pi/4}=|z|e^{\bi\pi/4}$. Therefore $\mb{Re}(\lam z)%=\mb{Im}(\lam z)
=\frac{\sqrt2}{2}|z|$.  Because
$$
|\mb{Re}(\lam b w)|\leq |b w|<|b||z|'\leq \sqrt2|b||z|,
$$
we deduce by \eqref{e:5.6d} that
$$
\mb{Re}(\lam \~w)\geq a\mb{Re}(\lam z)-|\mb{Re}(\lam b w)|\geq\sqrt2 \(a/{2}-|b|\)|z|>0.
$$
Similarly $\mb{Im}(\lam \~w)>0$. Hence $\lam \~w\in\o\bbC_+$. This and \eqref{e:5.6c} show that $\lam Au\in\o P$.

\subsection{Positivity of solutions of ODEs with  positive operators}
\vs
%Let us finally recall a fundamental  result concerning positivity of solutions of ODEs to conclude this section.
Let $E$ be a (real or complex) Banach space with a cone $P$. An operator  $ A\in\sL( E)$ is called {\em positive}, if $AP\subset P$.

Given a positive operator  $A\in\sL(E)$ and $\a\in\R$, consider in $E$  the system
$$
\dot x= A x+\a x, \hs x(0)=x_0.
$$
 Denote by $x(t;x_0)$ the solution of the system.

\bl\label{l:2.1}  If $x_0\in P$, then $x(t;x_0)\in   P$ { for all } $t\geq 0.$
\el
{\bf Proof.} This is a well known basic  fact. One  may consult \cite[pp. 60-61]{Henry} (Exercises 6-8) and \cite{Kry} for more general results.
   We include a  simple proof just for completeness and the reader's convenience.

Let $x_0\in  P$, and $t\geq0$. Since $A$ is a bounded linear operator, we have   %$e^{t  A}$ is the analytic semigroup generated by $ A$,
$$\ba{ll}
x(t;x_0)=e^{\a t}e^{t A }x_0,\hs\mb{where }\,e^{t A}=\lim_{m\ra\8}\sum_{k=0}^m \frac{t^k}{k!} A^k.\ea
$$
 The positivity of $A$ implies that   $ A^k x_0\in  P$ for all $k\geq 0$. Hence $\sum_{k=0}^m \frac{t^k}{k!} A^k x_0\in  P$ for $m\geq 1$. It follows by the  closedness of $ P$ that $e^{t A}x_0\in  P$. Therefore  $x(t;x_0)=e^{\a t} e^{t A}x_0\in  P$.  $\bx$

\section{Two Technical Lemmas}\label{s:3}

In this section we give two technical lemmas concerning a projective property of cones and the asymptotic behaviour  of linear ODEs, each of which may be of independent interest.

\subsection{Projective properness of cones}

\vs
Let $E$ be a Banach space (real or complex), and  $P$  a cone in $E$. One trivially verifies that if $E_0$ is a closed  subspace of $E$, then  $P_0:=P\cap E_0$ is a cone in $E_0$. % as long as $P_0\ne\{0\}$.

 \bd  Let $\cD: E=\oplus_{i=1}^n E_i$ be a  given direct sum decomposition of $E$.
  A cone  $ P$  is called  projectively proper (with respect to $\cD$), if  \be\label{e:2.0}\ba{ll}(\Pi_{i} P)\Cap\, (-\Pi_{i} P)=\{0\},\Hs 1\leq i\leq n,\ea\ee
  where $\Pi_i$ denotes the projection from $E$ to $E_i$.
\ed

By definition it is clear that  the cone $P=\{0\}$ is trivially projectively proper with respect to any direct sum decomposition $\cD$ of $E$.

Intuitively by ``projectively proper'' it means that {\em $P$ looks like a cone from any direction given by $\cD$}.  Generally a cone    may not possess this property itself. This can be seen from the following simple example.
 \Vs\noindent{\bf Example 3.1.}
Let $E=\R^2$, and let $E_1$ and $E_2$ be the $x$-axis and $y$-axis, respectively. Then $\cD: E=E_1\oplus E_2$ is a direct sum decomposition of $E$.
Set $P=\{(x,y)^T\in E:\,\,y\geq |x|\}$ (where $(x,y)^T$ denotes the transpose of $(x,y)$). Clearly  $P$ is a cone. However, since $\Pi_1 P=E_1$, $P$  is not projectively proper.
\Vs

Fortunately,  we have the following fundamental fact.

\bp\label{p:2.8}Let $P$ be a cone in $E$, $P\ne\{0\}$. Then for any direct sum decomposition  $\cD: E=\oplus_{i=1}^n E_i$ of $E$,   there is a nonempty  index set $J_0\subset \{1,2,\cdots,n\}$  such that
\be\label{e:pp}P_0:= P\cap   E_0\ne\{0\},\ee where   $ E_0=\oplus_{s\in J_0} E_{s}$; furthermore, $P_0$ is  a projectively proper cone  in $E_0$ (with respect to   $\cD_0:E_0=\oplus_{s\in J_0} E_{s}$\,).
\ep
{\bf Proof.} Let us argue by induction on $n$.
First, if $n=1$ then  the conclusion trivially holds true.
Suppose now that the conclusion holds true  with $n=N\geq 1$. We show  that it remains true  for $n=N+1$.

Thus we  let  $E=\oplus_{i=1}^{N+1} E_i$ be  a  direct sum decomposition of $E$.
If \eqref{e:2.0} holds true %for $(\Pi_{i} P)\Cap\, (-\Pi_{i} P)=\{0\}$
for all $i$, $1\leq i\leq N+1$, then $E_0:=E$ fulfills all the requirements with $\cD_0: E=\oplus_{i=1}^{N+1} E_i$. Hence we assume, without loss of generality, that
$$\ba{ll}\(\Pi_{  {N+1}} P\)\Cap\,\(-\Pi_{  {N+1}} P\)\ne\{0\}.\ea$$
Then  there is  an element  $x\ne0$ such that $\pm x \in\Pi_{{N+1}} P$. Let $\pm x=\Pi_{N+1} u^\pm$, where  $u^\pm\in P$.
Obviously $u^\pm\ne0$. We may write
$$u^+=y_1+x,\hs u^-=y_2-x,$$ where
 $y_1,y_2\in V_1=\oplus_{i=1}^{N} E_{i}$.
 Then $y=y_1+y_2\in V_1$. We observe that
\be\label{e:2.a}
y=(y_1+x)+(y_2-x)=u^++u^-\in  P.
\ee
We claim that $y\ne0$. Indeed, if  $y=0$ then  $y_2=-y_1$, and hence
$$
u^-=y_2-x=-y_1-x=-u^+\in  P.
$$
Since  $u^+\in P$, this leads to a contradiction.

We infer from the above argument that  $y\in P\cap V_1:=P_1$. Hence $P_1$ is a cone in $V_1=\oplus_{i=1}^{N} E_{i}$ with $P_1\ne\{0\}$. Now by the  induction hypothesis,  there is
a nonempty  index set $J_0\subset \{1,2,\cdots,N\}$  such that $P_0:= P\cap   E_0\ne\{0\}$ is a cone in $ E_0=\oplus_{s\in J_0} E_{s}$ which  is projectively proper. % with respect to the decomposition $\cD_0:E_0=\oplus_{i\in J_0} E_{i}$).
 $\bx$

%\Vs
%In what follows we give a simple example to help the reader have a better understanding to  the above proposition.
%
\br {\em In Proposition \ref{p:2.8} one may expect  that $J_0$ consists of exactly one index $j\in \{1,2,\cdots,n\}$.
  %there is an  index $j\in J$ such that $P_0:=P\cap E_j\ne\{0\}$,  in which case  the proposition then  holds true with $J_0=\{j\}$.
Unfortunately this fails to be true in  the general case.  For instance, let $E=\R^2$, and let $E_1$,  $E_2$ and $\cD$ be the same as in Example 3.1. Then the cone  $$\ba{ll}P=\{(x,y):\,\,x,y\geq 0,\,\,\frac{1}{2}x\leq y\leq x\}\ea$$  is projectively proper  with respect to $\cD$. However,  one trivially has   $P\cap E_j=\{0\}$   for both $j=1,2$. Therefore we see that $J_0=\{1,2\}$ is the only index set such that Proposition \ref{p:2.8} holds  for this example.}
\er

\subsection{Asymptotic behaviour of linear ODEs}
\vs
Let  $E$ be  a finite-dimensional complex Banach space, and  $ A\in\sL( E)$.
By basic knowledge in linear algebra,  $ E$ has a  unique {\em primary decomposition}
\be\label{e:dx1}\cD:\,\,\, { E}=\oplus_{i=1}^m E_i,\ee where  each  $E_i$ is a cyclic subspace of $ A$
generated by a generalized eigenvector  $\xi_i$ of $A$ corresponding to  an eigenvalue $\mu_i\in \sig(A)$,
$$E_i=\mb{span}\,\{\xi_i,\,( A-\mu_i) \xi_i\,,\cdots,\,( A-\mu_i)^{\kappa_i-1}\xi_i\},$$ where $\kappa_i=\mb{rank}(\xi_i)$ is the rank of $\xi_i$.
(It may happen that for $i\ne j$,  $\xi_i$ and $\xi_j$ correspond to the same eigenvalue, i.e., $\mu_i=\mu_j$.) We write
$$\mu_i=\a_i+\bi\b_i,\Hs i=1,2,\cdots, m.$$

%\subsection{Asymptotic behaviour of liner ODEs}
Denote by $x(t)=x(t;x_0)$ the solution of the initial value problem
\be\label{e:2.3}
\dot x= A x,\hs x(0)=x_0.\ee
 Then $x(t)=e^{t A }x_0$.
Thus   if $x_0\in \mb{GE}_{\mu}( A)$ with  $\mb{rank}(x_0)=\nu$, we have
\be\label{sode}
\ba{ll}x(t)&=e^{\mu t}e^{t( A-\mu) }x_0=e^{\mu t} \sum_{k=0}^\8 \frac{t^k}{k!}( A-\mu)^k x_0\\[1ex]
&=e^{\mu t}\(I+\frac{t}{1!}( A-\mu)+\cdots+\frac{t^{{\nu}-1}}{({\nu}-1)!}( A-\mu)^{({\nu}-1)}\)x_0.
\ea
\ee
\newpage
Let $J=\{1,2,\cdots,m\}$,
% For convenience, given $x\in E$, $x\ne0$, we set
%\be\label{e:3.5}J_\cD(x)=\{i\in J:\,\,\Pi_i x\ne 0\},\ee
and denote by   $\Pi_i$ the projection from $E$ to $E_i$. %Hence $x$ can be uniquely written as
%\be\label{e:3.6}\ba{ll}x=\sum_{i\in J_\cD(x)}u_i,\hs\mb{where }u_i=\Pi_i x\in E_i\sm\{0\}.\ea\ee

\bp\label{p:3.1}Let $x_0\in E$, $x_0\ne 0$.  We write $u_i=\Pi_i x_0$ \,$(i\in J)$, and let
$$J_*:=\{i\in J:\,\,u_i\ne 0\}.$$ Set  $\a=\max_{\,i\in J_*}\a_{i}$, and $J_\a=\{i\in J_*:\,\,\a_i=\a\}$. Put
$$\ba{ll}\nu=\max_{\,i\in J_\a}\nu_i,\hs\mb{where }\,\nu_{i}=\mb{\em rank}(u_{i}).\ea$$
Define
$J_0=\{i\in J_\a:\,\,\nu_{i}=\nu\}.$
 Then
$$\lim_{t\ra \8}\|t^{-(\nu-1)}e^{-\a t}{x(t)}-\gam(t)\|= 0,$$ where
\be\label{e:3.13}\ba{ll}\gam(t)=\sum_{i\in J_0} e^{\mathbf{i} \b_{i} t}w_{i},\hs w_i=\frac{1}{(\nu-1)!} ( A-\mu_i)^{\nu-1}u_i.\ea\ee

\ep

\br\label{r:w}{\em For each $i\in J_0$, since $u_i\in E_i\subset  \mb{GE}_{\mu_i}(A)$ and $\mb{ rank}(u_i)=\nu$, the vector $w_i$ given in \eqref{e:3.13} is an eigenvector of $ A$ pertaining to $\mu_i$. Note also  that  $w_i\in E_i$. Therefore  $w_i \,(i\in J_0)$ are linearly independent.}
\er

\noindent{\bf Proof of Proposition \ref{p:3.1}.}
 Set
$$J_1=J_\a\sm J_0,\hs J_2=J_*\sm J_\a.$$
Then  $J_0\cup J_1\cup J_2=J_*$. Hence we can rewrite $x_0$ as
$$\ba{ll}x_0=y_0+y_1+y_2,\hs\mb{where } y_j=\sum_{i\in J_j}u_{i}.\ea$$ Let $y_j(t)=e^{tA}y_j$. Obviously  $x(t)=e^{tA}x_0=y_0(t)+y_1(t)+y_2(t)$.

\vs
By \eqref{sode} we deduce that
\be\label{e:3.2b}
y_j(t)=\sum_{i\in J_j}e^{\a_i t}e^{\bi\b_i t}\(I+\frac{t}{1!}( A-\mu_i)+\cdots+\frac{t^{\nu_i-1}}{(\nu_i-1)!}( A-\mu_i)^{(\nu_i-1)}\)u_i
\ee
for $j=0,1,2$.
 %$$y_0(t)=\sum_{i\in J_0} e^{\a t}e^{\bi\b_i t}\(I+\frac{t}{1!}( A-\mu_i)+\cdots+\frac{t^{\nu-1}}{(\nu-1)!}( A-\mu_i)^{(\nu-1)}\)u_i.$$
Since $\a_i=\a$ and   $\nu_i<\nu$ for $i\in J_1$, and $\a_i<\a$ for $i\in J_2$, by \eqref{e:3.2b}  one trivially deduces that  $\left\|t^{-(\nu-1)}e^{-\a t}{y_j(t)}\right\|\ra 0$ as $t\ra\8$ for $j=1,2$.
%$$\left\|t^{-(\nu-1)}e^{-\a t}{y_j(t)}\right\|\leq M_jt^{-1},\Hs t\geq 1,\,\,j=1,2.$$
Now we observe that
$$y_0(t)=\sum_{i\in J_0} e^{\a t}e^{\bi\b_i t}\(I+\frac{t}{1!}( A-\mu_i)+\cdots+\frac{t^{\nu-1}}{(\nu-1)!}( A-\mu_i)^{(\nu-1)}\)u_i,$$
from which it can be easily seen that $t^{-(\nu-1)}e^{-\a t}{y_0(t)}=\gam(t)+R_\nu(t),
$
where  $\gam(t)$ is the function given in \eqref{e:3.13}, $R_\nu(t)\equiv 0$ if  $\nu=1$, and
$$
R_\nu(t)=\sum_{i\in J_0}\(e^{\bi\b_i t}\,t^{-(\nu-1)}\sum_{k=0}^{\nu-2}\frac{t^{k}}{k!}( A-\mu_i)^{k}u_i\)\,\,\,\,\mb{if }\nu\geq 2.$$  In either  case it is easy to see that   $R_\nu(t)\ra 0$ as $t\ra\8$.
Combining the above facts together one  immediately concludes the validity of the lemma. $\bx$

\section{The Finite-dimensional  KR theorem: A  General  Version and its Elementary Dynamical Proof via Linear ODEs}\label{s:4}

%The   main purpose  in this section is to prove via a completely self-contained dynamical approach  a general complex KR type theorem. As we have seen in Proposition \ref{p:2.7},  although tangential positivity is equivalent to resolvent positivity, we prefer to  state our results  in terms of tangentially positive operators. One consideration  in  doing so  is that,  in this present work the resolvent positivity of an operator is merely   used to guarantee the invariance of a cone $P$ under the corresponding semigroup. For tangentially positive operators,  this invariance property is naturally  guaranteed by the well-known viability theorem even in the case where $P$ is only a closed set. Therefore for such operators, a combination of  the viability theory and the dynamical approach  here might  allow us to extend the PF theory to more general  cases where cones  may  be replaced by other geometric objects.

%Now let us state and prove the main result in this section.
%Let $E$ be an $n$-dimensional {complex} Banach space, and $ P$ a cone in $E$. Let $ A\in\sL(E)$. Set $r=r(A)$.  For convenience in statement, if $r\in\sig(A)$ then   we  call $r$ and the corresponding eigenvectors the {\em principal eigenvalue}  and   {\em principal eigenvectors} of $A$, respectively.

%\subsection{A general complex Perron type theorem}

%A major part of this section is to prove the following generalized version of the  finite-dimensional Perron type   theorem  for  complex operators.

 Let $E$ be a finite-dimensional {complex} Banach space.  In this section we prove  a general  KR type theorem for positive  operators in $E$ via a completely self-contained dynamical approach based on linear ODEs.

Given an operator  $ A\in\sL(E)$, if $r:=r(A)\in\sig(A)$ then   we will  simply  call $r$ and the corresponding eigenvectors the {\em principal eigenvalue}  and   {\em principal eigenvectors} of $A$, respectively.

  \bt\label{t:4.1}%$($Complex PF theorem$)$
      Let  $ P$ be a  full cone in $E$, and  $ A\in\sL(E)$  a positive operator (with respect to $P$).   Then the following assertions hold:
        \benu
        \item[$(1)$] $r\in\sig(A)$ with  a principal  eigenvector  $w\in  P$.
        \item[$(2)$]  Let $\mu\in\sig(A)$. If \,$\mb{\em GE}_\mu(A)\cap P\ne\{0\}$, then $\mu\in[0,r]$ with a corresponding  eigenvector $v\in  P$.
            \item[$(3)$] No eigenvectors pertaining to other  eigenvalues $\mu\ne r$ are contained in $\o P$.
            %All eigenvectors of  eigenvalues $\mu$ with $\mu\ne r$ are contained in $ E\sm\o P$.
            \item[$(4)$] If $\o P\ne\emp$ and  contains a principal  eigenvector $w$,  the algebraic and the geometric multiplicities  of $r$ coincide.
            \eenu
\et
 \br {\em Assertion $(1)$ in the above theorem is a finite-dimensional complex version of the main part of the KR theorem (see e.g. \cite{Krein} or \cite[Theorem 1.1]{Du}). Our main contribution here is that we have given a different simple dynamical proof.  Assertion $(3)$  generalizes  a classical result on nonnegative matrices in \cite[Proposition 2.3 in Chap. I]{Schab} (see the corollary of the proposition) which states that for a nonnegative $n\X n$ matrix $A$, if $Aw=\lam w$ for some vector $w\in\R^n$ whose components   are strictly positive, then $\lam=r(A)$.

 To the best of our knowledge, assertions $(2)$ and $(4)$ are new  under the hypotheses in the theorem.}
 \er
\subsection{A weaker version of the  KR theorem}\vs
%In this section we prove the following theorem, asserting that every positive operator  has  real eigenvalues.
To prove Theorem \ref{t:4.1}, we start with a weaker version of the  KR theorem  without any additional assumptions on   $ P$ except that $P\ne\{0\}$.
\bt\label{t:4.3} Assume that $P\ne\{0\}$. Then  every  positive operator $A\in\sL(E)$ has a real  eigenvalue $\mu$ with a  corresponding eigenvector $w\in  P$.
\et
{\bf Proof.}
% We may assume $  P$ is total. Otherwise one can replace $E$ by $\mathbb{Y}= \ol{( P- P)}$.
Let   $\cD:E=\oplus_{i=1}^mE_i$ be the primary  decomposition of $E$ given in \eqref{e:dx1} with  each subspace $E_i$ corresponding to an eigenvalue $\mu_i=\a_i+\bi\b_i$.
By Proposition \ref{p:2.8} there is a nonempty index set  $J_0\subset \{1,2,\cdots,m\}$ such that $ P_0:=  P\cap   E_0\ne\{0\}$ is a projectively proper cone in $ E_0=\oplus_{s\in J_0} E_{s}$, i.e.,
\be\label{e:3.0}\ba{ll}
(\Pi_{s} P_0)\Cap \,(-\Pi_{s} P_0) =\{0\},\Hs s\in J_0,\ea\ee
where $\Pi_i$ denotes the projection from $E$ to $E_i$.

 Lemma \ref{l:2.1} asserts that $e^{tA}P\subset P$ for $t\geq 0$.  Since $E_0$ is invariant under the semigroup $e^{tA}$ (i.e. $e^{tA}E_0\subset E_0$), we trivially have
\be\label{e:3.2}
 e^{tA}  P_0= e^{tA}(  P\cap  E_0)=( e^{tA}  P)\cap ( e^{tA} E_0)\subset   P\cap  E_0=  P_0.
\ee
 Pick an $x_0\in   P_0$, $x_0\ne0$, and let $x(t)=e^{tA}x_0$.  Thanks to Proposition \ref{p:3.1} and  Remark \ref{r:w},  there exist $\a,\nu\in\R$ with $\nu\geq 1$ and a nonempty index set $J_\a\subset J_0$ such that
\be\label{e:3.1}
\lim_{t\ra \8}||t^{-(\nu-1)}e^{-\a t}x(t)-\gamma(t)||=0,
\ee
where  $\gamma(t)=\sum_{s\in J_\a}  e^{\mathbf{i}\b_{s} t}w_{s},$
   and $w_{s}\in E_s$ is an  eigenvector of $A$ corresponding to $\mu_s$ for each $s\in J_\a$. Furthermore, $\mu_s$ ($s\in J_\a$) share the same real part $\a$, that is,   $\mu_{s}=\a_s+\mathbf{i}\b_{s}=\a+\mathbf{i}\b_{s}$.

\vs

We claim   that \be\label{e:gam}\gam(t)\in  P_0,\Hs \A\,t\in\R.\ee
%For this purpose, let $y_s(t)=\Pi_{s} y(t),$ where $y(t)=t^{-(\nu-1)}e^{-\a t}x(t).$
Indeed, by \eqref{e:3.2} we have  $x(t)\in   P_0$ for $t\geq 0$. Therefore
%\be\label{e:3.2}y(t):=t^{-(\nu-1)}e^{-\a t}x(t)\in   P_0,\Hs t\geq 0.\ee
$y(t):=t^{-(\nu-1)}e^{-\a t}x(t)\in   P_0$ for $t> 0$.
Noticing that  $\gam(t)$ is an almost periodic function in the terminology in  Corduneanu \cite{Cond},
for every  $t\in\R$ fixed, by the basic properties of almost periodic functions (see Property B in \cite[pp. 24]{Cond}) it is easy to deduce that there is a sequence $t_n\ra \8$ such that
$$
\lim_{n\ra\8}\|\gam(t_n)-\gam(t)\|=0.
$$
Thus  by \eqref{e:3.1} we have  $\|y(t_n)-\gam(t)\|=0$ as $n\ra\8$. Since $y(t_n)\in P_0$ for $t_n>0$,  the closedness of $P_0$   implies   that $\gam(t)\in P_0$. Hence the claim holds true.

By \eqref{e:gam} one concludes that for each $s\in J_\a$,
\be\label{e:gam1}
\gam_s(t):=\Pi_s\gam(t)=e^{\mathbf{i}\b_{s} t}w_{s}\in \Pi_s P_0,\Hs \A\,t\in\R.
\ee

Now we show that  $\b_{s}=0$ for all $s\in J_\a$, and hence $\sig_0:=\{\mu_s:\,\,s\in J_\a\}$ consists of exactly one real eigenvalue $\mu=\a$.
We argue by contradiction and suppose  that   $\b_{s}\ne 0$ for some $s\in J_\a$. Pick  a $\tau>0$ such that
$
e^{\mathbf{i}\b_{s} \tau}= -1.
$
Then
$$
\gam_s(\tau)= e^{\mathbf{i}\b_{s} \tau}w_{s}=- w_{s}=-\gam_s(0).
$$ This and   \eqref{e:gam1} assert   that $\pm \gam_s(0)\in  \Pi_s P_0$, which contradicts \eqref{e:3.0}.

Let us finally verify that $A$  has an eigenvector $w\in P_0\subset P$. Since $\sig_0=\{\a\}$,  all the eigenvectors $w_s$ ($s\in J_\a$) correspond to the same eigenvalue $\mu=\a$. Consequently
$w:=\sum_{s\in J_\a} w_{s}$ is an eigenvector of $A$ pertaining to $\a$. Note that
$$\gam(t)=\sum_{s\in J_\a}   e^{\mathbf{i}\b_{s} t}w_{s}\equiv w.$$ Because  $t^{-(\nu-1)}e^{-\a t}x(t)\in P_0$ for all $t\geq 0$, by \eqref{e:3.1} we deduce that $w\in P_0$. This completes the proof of the theorem. $\bx$

\br {\em It is interesting to mention that in the proof of Theorem \ref{t:4.3} we have only used the fact that the semigroup $e^{tA}$ is positive, meaning that
$e^{tA}P\subset P$ for all $ t\geq 0.$
(The positivity of the operator $A$  is just used to guarantee the positivity of $e^{tA}$.) We remark  that the generator of  a positive semigroup may not be positive itself, and therefore Theorem \ref{t:4.3} actually holds true for a wider class of operators.
A somewhat trivial  example is the  real square matrix $A=(a_{ij})_{n\X n}$ with $a_{ij}\geq 0$ for $i\ne j$. Obviously $A$ is not positive if $a_{ii}<0$ for some $i$. However, if we take a number $\a>\max_{1\leq i\leq n}|a_{ii}|$, then the matrix $A_\a:=A+\a I$ has nonnegative entries and hence preserves the cone  $P:=\R^n_+$ in $E=\R^n$. It then immediately follows from Lemma \ref{l:2.1} that $e^{tA}P\subset P$ for $t\geq 0$.

 It is also desirable to obtain analogues of Theorems \ref{t:4.1} and \ref{t:6.1} etc. for generators of positive semigroups  in infinite-dimensional Banach spaces. However, this needs some further careful investigations, which may deserve to be carried out in a future work.
 }
\er
\br {\em Some elegant descriptions concerning operators which generates cone-preserving semigroups have been given  in Kanigowski and Kryszewski \cite[Section 3]{Kry} in terms of resolvent invariance of  cones and a tangency condition on operators  with respect to the cones.   Based on these fundamental results  the authors proved some generalized KR theorems for the generator of a cone-preserving  semigroup in a Banach space  without assuming that the cone induces a Banach lattice structure into the  space.}
\er

%%%%%%%%%%%%%%%%%%%%%%%%%%%%%%

\subsection{Proof of Theorem \ref{t:4.1}}
%This section is devoted to  the  proof of a finite-dimensional version of the well-known KR theorem.
\vs
 We split  the proof of Theorem \ref{t:4.1} into two independent lemmas.

 \bl\label{l:4.4}Let  $P$ be a cone in $E$, and   $ A\in\sL(E)$  a positive operator with respect to $P$.
If $ P$ is full, then $r\in\sig(A)$ with a principal   eigenvector  $w\in  P$.
        \el

\vs\noindent{\bf Proof.}  If $r=0$ then  $\sig( A)=\{0\}$. By virtue of Theorem \ref{t:4.3}  the conclusion  trivially holds true. Thus we assume  $r>0$. It can be assumed  that $r=1$; otherwise one may   argue with  $ A_1=r^{-1} A$ in place of $ A$.
 %The basic idea is to  reduce the problem  to a simpler one  where $\sig(A)\subset \bbS_1$ and then apply Theorem \ref{t:4.3}.

Take a  $\de>0$ sufficiently small so that all the eigenvalues $\mu$ of $ A$ with $|\mu|>1-2\de$ are contained in  $\bbS_1$.
% see Fig. \ref{fg3-2}.
%\begin{figure}[h!]    \begin{center}       \includegraphics[width=4.5cm]{fig3-2.pdf}  \caption{Ditribution of eigenvalues}    \label{fg3-2}       \end{center}     \end{figure}
Set
$$
\sig_0=\{z\in \sig( A):\,\,|z|\leq 1-2\de\},\hs \sig_1=\{z\in \sig( A):\,\,|z|=1\}.
$$
 Obviously $\sig( A)=\sig_0\cup \sig_1$. Let $ E=  E_0\oplus   E_1$ be the  corresponding  decomposition of $E$. Then both $ E_0$ and $ E_1$ are invariant subspaces of $ A$. %(Such a decomposition can be easily obtained by using the primary one in \eqref{e:dx1}.)
Let $P_1=P\cap   E_1$. Then $P_1$ is a cone in $E_1$. We show that \be\label{e:3.0a}  P_1\ne\{0\} .\ee

Take a number $\rho$ with $1<\rho\leq(1-\de)/(1-2\de)$, and set % such that $\rho\sig_0\subset \{z\in\bbC:\,\,|z|\leq 1-\de\}.$ Set
$B=\rho  A$. Then
$$\sig(B|_{ E_0})=\rho\sig_0\subset \{z\in\bbC:\,\,|z|\leq 1-\de\},$$
$$
\sig(B|_{ E_1})=\rho\sig_1\subset \{z\in\bbC:\,\,|z|=\rho\}.
$$
 Hence we deduce that  \be\label{e:4.2}\mb{$||B^k x||\ra 0$} \,\,\,(x\in   E_0),\hs \mb{and }\,\mb{$||B^k x||\ra \8$}\,\,\,(x\in  E_1\sm\{0\}).\ee

 We may assume that $\mb{dim}(E_0)\geq 1$; otherwise $ E=E_1$ and \eqref{e:3.0a} readily holds. Since $ P$ is full, we deduce   that
 $ P\sm   E_0\ne \emp$. Pick an $x\in   P\sm   E_0$, $x\ne 0$. We may write  $x=x_0+x_1$, where $x_i\in   E_i$. Then  $x_1\ne0$.  Note that   $B^kx_i\in E_i$ for all $k$. By the first relation in \eqref{e:4.2} we have $B^k x_0\ra0$. Therefore
\be\label{e:3.9}
\lim_{k\ra\8}d(B^k x,  E_1)=\lim_{k\ra\8}d(B^k x_0+B^k x_1,  E_1)=0.
\ee
Now if $  P\cap   E_1=\{0\}$, then because  $B^kx=\rho^kA^kx\in P$ for $k\geq 1$, by Lemma \ref{l:2.4} we deduce  that $B^k x\ra0$. Consequently   $B^k x_1=B^k x-B^k x_0\ra0$. This leads to a contradiction. Hence  \eqref{e:3.0a} holds true.

\vs
 Let $A_1= A|_{E_1}$. Then $A_1$ is positive in $E_1$ (with respect to  $P_1$) with  $\sig(A_1)\subset \bbS_1$.
 By virtue of Theorem \ref{t:4.3}  one  concludes  that $\sig(A_1)$ contains a real eigenvalue $\mu$ with  a corresponding  eigenvector $w\in P_1\subset P$.  Since $\mu\in\bbS_1$, we necessarily have $\mu\in\{1,-1\}$.

We claim that $\mu\ne-1$, and consequently $ \mu=1$, which finishes the proof of the lemma.
Indeed, suppose  $\mu=-1$. Then since $w\in P$, by the positivity of $A$ we would have  $-w=A_1 w=Aw\in P$. This leads to a contradiction.
 $\bx$

%By the definition of $\a$, this also implies that $\a=\mb{Re}\,\lam=1$.  Since  $\Sigma_\a\subset \sig(B)\subset \bbS_1$ and $\mb{Re}\,\mu=\a$ for $\mu\in\Sigma_\a$, it can be  easily seen  that $\Sigma_\a=\{\lam\}$. $\bx$

 \bl\label{l:4.5}Let  $P$ be a cone in $E$, and   $ A\in\sL(E)$  a positive operator with respect to $P$.
    If  $r\in\sig(A)$ with a principal  eigenvector  $w\in  P$, then the following assertions hold:
        \benu
               \item[$(1)$] Given $\mu\in\sig(A)$, if\,  $\mb{\em GE}_\mu(A)\cap P\ne\{0\}$, %$ P$ contains a generalized eigenvector $\zeta$ of $A$ pertaining to $\mu$,
                   then $\mu\in[0,r]$ with a   corresponding eigenvector $v\in  P$.
            \item[$(2)$] The eigenvectors of any other   eigenvalue $\mu\ne r$ are contained in $ E\sm\o P$.
            \item[$(3)$] If $\o P\ne\emp$ and contains a principal eigenvector, then the algebraic and the geometric multiplicities  of $r$ coincide.
            \eenu
\el
{\bf Proof.}  (1) We may assume  $\mu\ne r$. %Assume that $ P$ contains a generalized  eigenvector $\zeta$ corresponding to $\mu\in\sig(A)$. We show  that  $\mu\in[0,r)$ with an eigenvector $\eta\in  P$.
     Let $\mu=\a+\mathbf{i}\b$. Pick a $\zeta\in \mb{GE}_\mu(A)\cap P$, $\zeta\ne0$. Denote by $x(t)$ the solution of the system  $\dot x= A x$ with $x(0)=\zeta$. Then
$$\ba{ll}
x(t)=e^{(\a+\mathbf{i}\b) t}\(I+\frac{t}{1!}( A-\mu)+\frac{t^2}{2!}( A-\mu)^2+\cdots +\frac{t^{\nu-1}}{(\nu-1)!}( A-\mu)^{\nu-1}\)\zeta,\ea
$$
 where $\nu=\mb{rank}(\zeta)$. We observe that
 $$
\|t^{-(\nu-1)}e^{-\a t}x(t)-\gam(t)\|\ra 0\hs\mb{as }\,t\ra\8,
$$
where $\gam(t)=e^{\mathbf{i}\b t}\eta$, and $\eta=\frac{1}{(\nu-1)!}( A-\mu)^{\nu-1}\zeta$ is an eigenvector of $ A$ associated with  $\mu$. Since $x(t)\in  P $ (by Lemma \ref{l:2.1}), as in the proof of Theorem  \ref{t:4.3} it can be shown that $\gamma(t)\in  P $ for all $t\in \R$.
In particular,
\be\label{e:4.5}\eta=\gam(0)\in  P .\ee

\vs We show   that $\b= 0$, therefore  $\mu=\a\in\R$.
Indeed, suppose  $\b\ne 0$. Then one can find a   $\tau>0$ such that $e^{\mathbf{i}\b \tau}=-1$. Thus $-\eta=e^{\mathbf{i}\b \tau}\eta =\gam(\tau)\in   P $, which  contradicts \eqref{e:4.5}.

To prove assertion (1), it remains to check  that $\mu\geq 0$. First, by \eqref{e:4.5} and the positivity of $ A$,  $\mu\eta= A\eta \in  P$. Now if $\mu<0$, then
$
-\eta=\frac{\mu}{|\mu|}\eta=\frac{1}{|\mu|}\(\mu\eta\)\in  P ,
$
which yields   a contradiction. Thus we see that   $\mu\geq 0$.

\vs
(2) \,Let $v\in  P$ be an eigenvector corresponding  to an eigenvalue $\mu\ne r$. We need to check  that $v\not\in\o P$.

 We argue by contradiction and suppose that $v\in\o P$ (hence $\o  P\ne\emp$).
Then by  assertion (1) we deduce that   $\mu\in[0,r)$. Let $w$ be the principal eigenvector  given by assumption, and let   $\pi$ be the    {\em real plane }  spanned by $w$ and $v$,
$$\pi=\{a w+b v:\,\,a,b\in\R\}.$$ Obviously $\pi$  is  isomorphic to $\R^2$.  Note that the cone    $ P_\pi:=P \cap  \pi\ne\{0\}$. Denote by $\o{P}_\pi$ the interior of $  P _\pi$ in $\pi$. One easily sees that $\o{P}_\pi\ne\emp$.

As  $w$ and $v$ are eigenvectors of $ A$ associated with real eigenvalues, we trivially  have    $ A \pi\subset \pi$. Let $ A_\pi=A|_\pi$.
Then  $  A_\pi$ is  a {\em real\,} linear operator on $\pi$. Fix a real number $\lam$ with $\mu<\lam<r $ and consider the planar system
\be\label{e:3.7}
\dot y=(  A_\pi-\lam)y,\hs y(0)=y_0.
\ee
By Lemma \ref{l:2.1}  the cone  $  P _\pi$ is invariant under system \eqref{e:3.7}.

\begin{figure}[h!]
    \begin{center}
       \includegraphics[width=6cm]{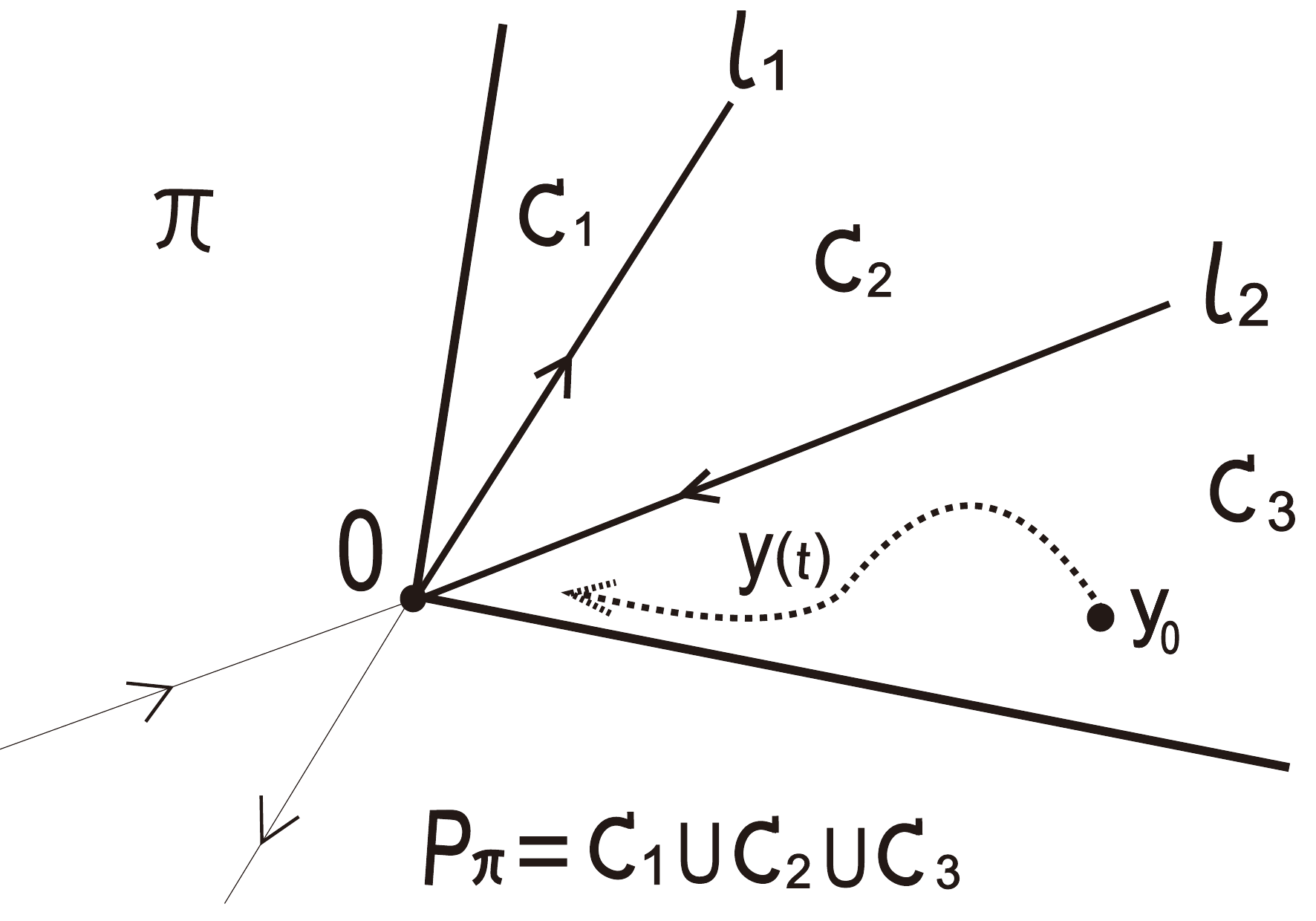}\vs  \caption{$C_2$ and $C_3$ are solid cones in $\pi$}
    \label{fg3-1}
    \end{center}
     \end{figure}

The two lines $l_1:=\{s w:\,\,s\in\R\}$ and $l_2:=\{s v:\,\,s\in\R\}$ split $  P _\pi$ into three subcones $C_1,C_2$ and $C_3$, where $C_2$ denotes the cone with boundary $\pa C_2=\(l_1\cup l_2\)\cap P_\pi$; see Fig. \ref{fg3-1}.
Since $l_1\ne l_2$,  $C_2$ is a solid cone in $P_\pi$.  Because   $v\in\o{ P} _\pi$, it is also easy to   see that  $\o C_3:=\mb{int}_\pi C_3\ne\emp$.
  %Since $l_1\ne l_2$, it is also clear that $C_2$ is solid. %It can be easily seen that $C_1\cap C_3=\{0\}.$

  The operator $  A_\pi-\lam$ has two eigenvalues $\mu_1:=r-\lam>0$ and $\mu_2:=\mu-\lam<0$ corresponding to   eigenvectors $w$ and $v$, respectively. Hence  $l_1$ and $l_2$ are unstable and stable manifolds of system \eqref{e:3.7}, respectively.
  Since any solution of \eqref{e:3.7} can not cross the lines $l_1$ and $l_2$,
  we infer from the invariance  of $P_\pi$ with respect to  \eqref{e:3.7}  that every  subcone $C_i$ ($i=1,2,3$) is preserved  under  \eqref{e:3.7}.

   Take  a $y_0\in \o C_3$. Write  $y_0=a w+b v$. Then $a\ne0\ne b$. The solution of \eqref{e:3.7} reads as
$
y(t)=ae^{\mu_1 t}w+b e^{\mu_2 t}v.
$
Since $b e^{\mu_2 t}v\ra0$  as $t\ra \8$, we  deduce that
\be\label{e:4.0}\lim_{t\ra\8}|y(t)|=\lim_{t\ra\8}|ae^{\mu_1 t}w|=\8,\ee
and \be\label{e:3.8}
\lim_{t\ra\8}d(y(t), l_1)=0.
\ee
On the other hand, the invariance of $C_i$ with respect to \eqref{e:3.7} implies that $y(t)\in C_3$ for all $t\geq 0$. Because   $l_1\cap C_3=\{0\}$,  by  \eqref{e:3.8} and Lemma \ref{l:2.4} we conclude  that $y(t)\ra 0$ as $t\ra\8$, which contradicts \eqref{e:4.0} and finishes the proof of (2).
\vs

(3) \,Assume that $\o P\ne\emp$ and   contains a principal eigenvector $\omega$. We need to check  that rank$(\xi) =1$ for every $\xi\in\mb{GE}_r(A)$, $\xi\ne0$.

We argue by   contradiction and suppose  that  $\nu:=\mb{rank}(\xi)\geq 2$ for some $\xi\in\mb{GE}_r(A)$, $\xi\ne0$. Let $\xi_0=( A-r)^{\nu-2}\xi$. Clearly  $\xi_0\in\mb{GE}_r(A)$ with $\mb{rank}(\xi_0)=2$.
Since  $\omega\in\o P$, there exists  $\de>0$ such that  $u:=\omega+z \xi_0\in  P$ for $z\in\bbC$ with $|z|< \de$. Note that $u\in \mb{GE}_r(A)$ with $\mb{rank}(u)=2$.
  The same argument as in the verification of \eqref{e:4.5} then applies to show that $$v:=(A-r)u=z ( A-r)\xi_0 \in  P.$$ Taking $z=\pm\frac{\de}{2}$, one concludes that $\pm\frac{\de}{2}( A-r)\xi_0\in  P$, which leads to a contradiction because $( A-r)\xi_0=( A-r)^{\nu-1}\xi\ne0$. $\bx$

\Vs
\noindent{\bf Proof of Theorem \ref{t:4.1}.} Combining Lemmas \ref{l:4.4} and \ref{l:4.5} together, one immediately concludes the validity of all assertions in the theorem.  $\bx$

\section{The Finite-dimensional KR theorem: A Strong Version}\label{s:5}
%In this section we state and prove a strong version of the  KR theorem,  extending  a major part of the classical KR theorem   to complex operators. A simple  dynamical scheme for  calculating  principal eigenpairs will also be presented.

%\subsection{A strong complex  version of the KR theorem}
\vs
Let $ E$ be a finite-dimensional complex Banach space, and $P$ a cone in $E$.
 Our main purpose  in this section is to prove   the following theorem.

\bt\label{t:5.2}Let   $ A\in\sL(E)$ be a  positive operator.
 Suppose $P$ is solid, and  that $ A$  is rotationally strongly positive. Then the following assertions hold:
\benu
\item[$(1)$] $r:=r(A)\in\sig(A)$ with a principal   eigenvector $w\in P$.
\item[$(2)$] The algebraic and the  geometric multiplicities  of  $r$ coincide; moreover, \be\label{e:5.16}\ba{ll}(\bbS_1 \xi)\cap \o P\ne\emp\ea\ee for any principal eigenvector $\xi$.
\item[$(3)$] $\mb{\em GE}_\mu(A)\cap P=\{0\}$ for  any other  eigenvalue   $\mu\ne r$.
\item[$(4)$] $|\mu|<r$ for all $\mu\in \sig(A)\sm\{r\}$.
\eenu
\et

\noindent{\bf Proof.} Assertion (1) directly follows from Theorem \ref{t:4.1} (1). So we only need to verify (2)-(4).

We claim that $r>0$. Indeed, suppose $r=0$. Let $w\in P$ be the principal eigenvector given in assertion (1). Then $Aw=0$. Consequently $\bbS_1(Aw)=\{0\}$. Hence $\bbS_1(Aw)\,\cap \o P=\emp$,  which contradicts the rotational strong positivity of $A$. Hence the claim holds true.

As before, we may assume that  $r=1$.

Let $w$ be the principal eigenvector in assertion (1).
By the  rotational strong positivity of $A$, there exists $\lam_1\in \bbS_1$ such that $\lam_1 A w\in \o{P}$. Set
$\xi_1=\lam_1 w$. Clearly $\xi_1$ is a principal  eigenvector. Note that
\be\label{e:5.2}\xi_1=A\xi_1=\lam_1 A w\in\o{P}.\ee

In what follows we first prove assertion  (4).  Let us  argue by way of contradiction and suppose that $ A$ has an eigenvalue $\mu\ne 1$ with $|\mu|=1$.  Let $\xi_2$ be an eigenvector of   $\mu$. Then  $z\xi_2$ is an eigenvector of    $\mu$ for every $z\in \bbS_1$. Set
$$C_1= \bbS_1\xi_1,\hs C_2=\bbS_1\xi_2.$$
We claim  that $C_2\cap {P}=\emp$. Indeed, if $z\xi_2\in {P}$ for some $z\in \bbS_1$, then since $\mu\ne 1$,  by Theorem \ref{t:4.1} (2) we deduce  that $\mu\in[0,1)$,  which leads to a contradiction.

Now that   $C_2\cap {P}=\emp$, by the compactness of $C_2$ we have
\be\label{e:3.0d}
\inf_{v\in C_2}d(v,{P}):=\de>0.
\ee

For each $t\geq 0$, set
$$
M(t):=C_1+tC_2=\{v_1+t v_2:\,\,\,v_i\in C_i\}.
$$
Then $M(t)$ is a compact subset of $E$ for each $t$. Furthermore, one  easily sees that $M(t)$ is continuous in $t$ in the sense of the Hausdorff distance $\de_H(\.,\.)$ given by
$$
\de_H(M,N)=\max\(d_H(M,N),\,d_H(N,M)\)
$$
for any bounded sets $M,N\subset E$, where $d_H(M,N)=\sup_{x\in M}d(x,N)$.
%, meaning that  for every $t\geq 0$ and  $\ve>0$, there is $\de>0$ such that for any $t'\geq 0$ with $|t'-t|<\de$, we have
%$$M(t')\subset \mB_\ve(M(t)),\hs\mb{and }\, M(t)\subset \mB_\ve(M(t')),$$
%where $\mB_\ve(M):=\{x\in E:\,\,d(x,M)<\ve\}$ is the $\ve$-neighborhood of a set $M\subset E$ in $E$.

Noticing that $t^{-1}M(t)=t^{-1}C_1+C_2$,  by \eqref{e:3.0d} it is easy to see that
if $t>0$ is sufficiently large, then $(t^{-1}M(t))\cap P=\emp$. Hence
 $$M(t)\cap P=t\((t^{-1}M(t))\cap (t^{-1}P)\)= t\((t^{-1}M(t))\cap P\)=\emp.$$
 On the other hand, since   $\xi_1\in\o{P}$, we see that if $t>0$ is sufficiently small then $\xi_1+tv_2\in\o{P}$ for all $v_2\in C_2$. Therefore  $M(t)\cap \o{P}\ne \emp$.
\vs
Define
$
\tau=\inf\{t>0:\,\,M(t)\cap P=\emp\}.
$
Clearly $\tau>0$. Since  $M(t)$ is continuous  in $t$, by compactness of $M(t)$  we deduce that $M(\tau)\cap P\ne\emp$.
We claim  that
\be\label{e:mt}M(\tau)\cap P \subset \pa P.\ee
 Indeed, if $M(\tau)\cap \o{P}\ne\emp$ then the continuity of $M(t)$ in $t$ implies that  $M(t)\cap \o{P}\ne\emp$ for $t\in [\tau,\tau+\ve]$ provided that  $\ve>0$ is sufficiently small. This contradicts the definition of $\tau$.

Take  a $u\in M(\tau)\cap P$. By the  rotational strong positivity of $A$, there is a $z'\in \bbS_1$ such that $z' A u\in\o{P}$. Let $u=v_1+\tau v_2$ for some $v_i\in C_i$, $i=1,2$. Noticing that   $Av_1=v_1$ and $Av_2=\mu v_2$, we find that
$$
z' A u =z' A (v_1+\tau v_2)=z'v_1+\tau (z'\mu) v_2.
$$
Since $|z'|=|z'\mu|=1$, we see that $z'v_1\in C_1$, and $(z'\mu) v_2\in C_2$. Therefore $z' A v\in M(\tau)$. But this and $z' A u\in\o{P}$ contradict \eqref{e:mt}.
\vs

Now we  turn to the verification  of   assertion (2). The conclusion that  $r$ has the same  algebraic and geometric multiplicity is a simple consequence of \eqref{e:5.2} and Theorem \ref{t:4.1} (4). So we only need to verify  \eqref{e:5.16}.

Let $\xi$ be a principal  eigenvector. %We claim   that
%\be\label{e:5.18}\ba{ll}(\bbS_1\xi)\cap P\ne\emp.\ea\ee Indeed,
If  $(\bbS_1\xi)\cap P=\emp$ then one can repeat  the same  argument as in the verification of   assertion (4) with $\xi_2$ therein replaced by $\xi$ and  obtain a contradiction. Therefore we deduce that $(\bbS_1\xi)\cap P\ne \emp$.
Thus there is  a $z\in \bbS_1$ such that $z\xi\in P$. Note that $w'=z\xi$ is a principal eigenvector. Now  repeating the same  argument as in leading to \eqref{e:5.2} with $w$ therein replaced by $w'$, we deduce  that
$\lam'w'\in \o P$ for some $\lam'\in\bbS_1$, i.e., $\eta\xi\in \o P$, where $\eta=z\lam'$. Since $\eta\in\bbS_1$, this finishes the proof of  \eqref{e:5.16}.
%by rotational strong positivity of $A$ we find that
%\be\label{e:5.19}\ba{ll}\emp\ne \bbS_1(A(z\xi))\Cap\o P=\bbS_1(zr\xi)\Cap\o P=r(z\bbS_1\xi)\Cap\o P=r(\bbS_1\xi)\Cap\o P.\ea\ee Since $r$ is a positive real number, \eqref{e:5.19} implies  that $(\bbS_1\xi)\cap \o P\ne\emp$.
\vs
It  remains to prove assertion (3).  We argue by contradiction and suppose  that $\mb{GE}_\mu(A)\cap P\ne\{0\}$ for some  $\mu\in\sig(A)\sm\{ r\}$. Then by Theorem \ref{t:4.1} (2), there is  a corresponding  eigenvector $\xi\in P$. By the rotational strong positivity of $A$ we deduce that $zA\xi\in \o P$ for some $z\in \bbS_1$. That is, $z\mu\xi\in\o P$. But this contradicts Theorem \ref{t:4.1} (3) because    $v:=z\mu\xi$ is an eigenvector of  $\mu$.  $\bx$
\br {\em Theorem \ref{t:5.2} can be seen as a complex version of the classical Perron theorem for positive matrices (or the finite-dimensional KR theorem). The  main contribution here  is that we have reformulated the Perron  theorem in a complex fashion by using the notion of rotational strong positivity instead of that of strong positivity widely used for real operators and given  a completely self-contained proof for it that is different from the existing ones in the literature.}
\er

\br {\em It remains an open question  whether  the principal eigenvalue  $r$ is geometrically simple  under the hypotheses of Theorem \ref{t:5.2}.}\er

 \section{The KR theorem: Infinite-dimensional Versions}\label{s:6}
\vs
We now state and prove some   infinite-dimensional   versions of Theorems \ref{t:4.1} and \ref{t:5.2} for complex operators.

Let $ E$ be an infinite-dimensional complex Banach space, and   $ P$ a cone in $ E$.
Let   $ A\in\sL(E)$. For notational simplicity, we write $r_e(A)=r_e$, and   $r(A)=r.$

\begin{figure}[htbp]
     \centering
     \subfigure[]{
     \begin{minipage}{5cm}
     \centering
     \includegraphics[width=5.5cm]{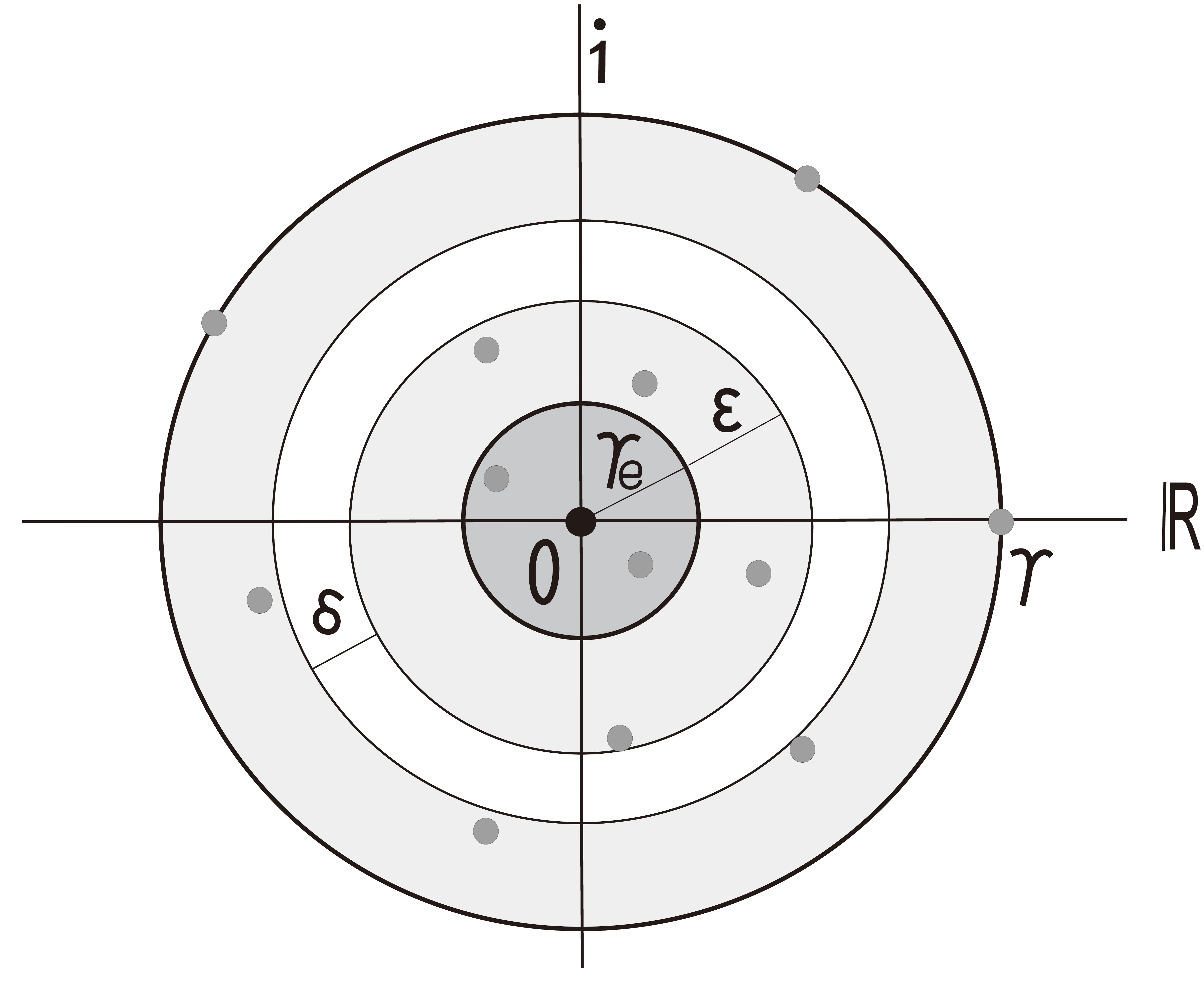}
     \end{minipage}
     }\Hs
     \subfigure[]{
     \begin{minipage}{5cm}
     \centering
     \includegraphics[width=5.5cm]{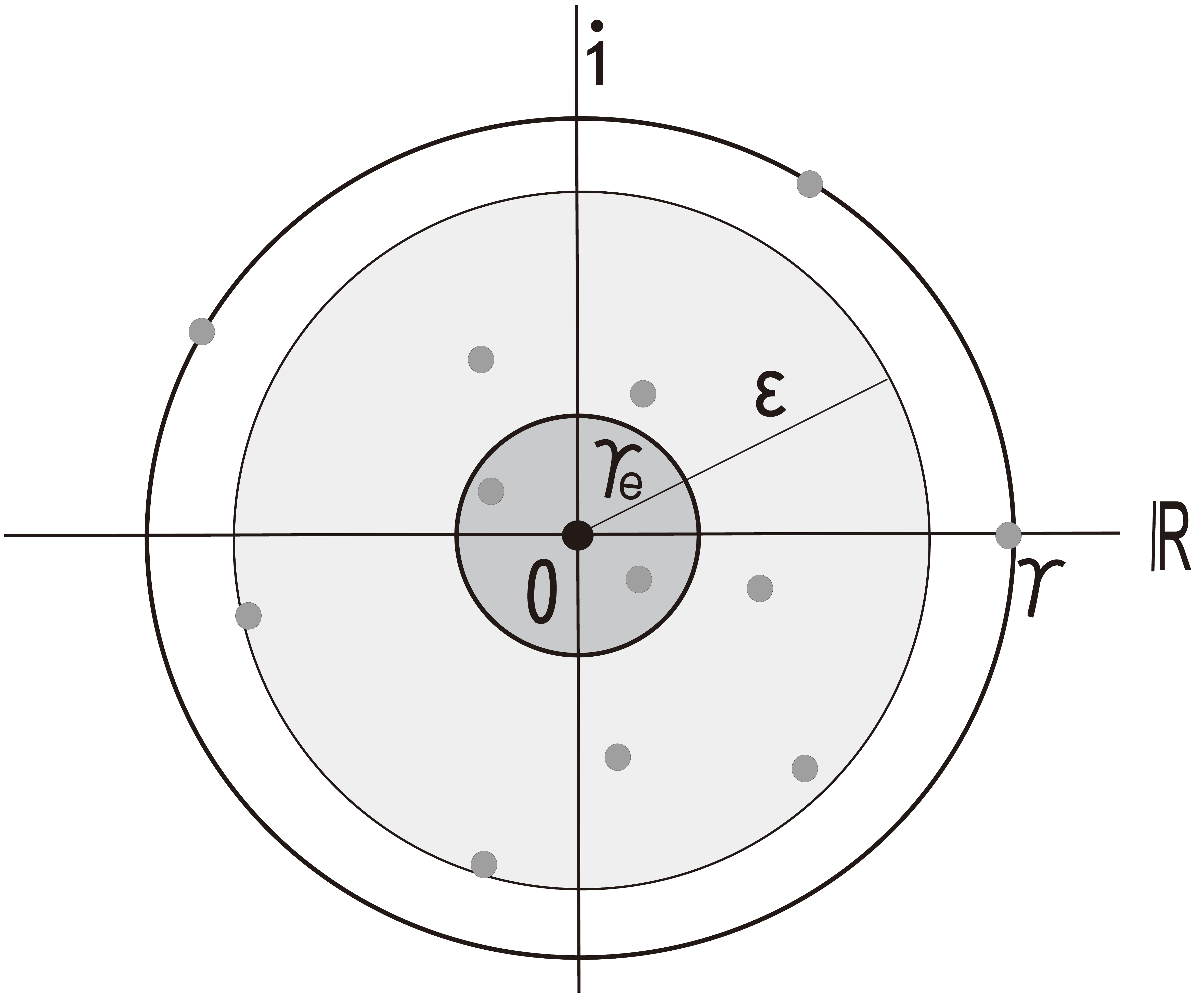}
     \end{minipage}
     }
     \caption{Distribution of the Spectrum}
     \label{fig:4.1}
     \end{figure}

   \bt\label{t:6.1} Let  $A$ be  a positive operator with $r_e<r$.
 Assume that  $ P$ is full.   Then  the following assertions hold true:
        \benu
        \item[$(1)$] $r$ is an eigenvalue of $ A$ with a principal eigenvector  $w\in  P$.
        \item[$(2)$] Let $\mu\in\sig(A)$, $|\mu|>r_e$. If $  P$ contains a generalized eigenvector pertaining to $\mu$, then $\mu\in[0,r]$ with  a corresponding  eigenvector $v\in  P$.
            \item[$(3)$] No eigenvectors of $ A$ corresponding to  other  eigenvalues  $\mu\ne r$ with   $|\mu|>r_e$ are contained in $ \o P$.
            \item[$(4)$] If $\o P\ne \emp$ and  contains a principal eigenvector, then the algebraic and geometric multiplicities  of $r$ coincide.
            \eenu
\et
{\bf Proof.} % Let us first verify assertion (1).
For any $\ve\in(0,\,r-r_e)$,  the region $\{z\in\bbC:\,\,r_e+\ve\leq |z|\leq r\}$ contains only a finite number of elements in $\sig( A)$ (see Section \ref{s:2.3}).  Hence  there is $\de>0$ such  that $\sig( A)\cap \{z\in\bbC:\,\,r_e+\ve<|z|< r_e+\ve+\de\}=\emp;$
see Fig. \ref{fig:4.1} (a). Set
\be\label{esig}
\ba{ll}
\sig_0(\ve)=\sig( A)\Cap \,\{z\in\bbC:\,\,|z|\leq r_e+\ve\},\\[1ex]
 \sig_1(\ve)=\sig( A)\Cap \,\{z\in\bbC:\,\,|z|\geq r_e+\ve+\de\}.\ea
\ee
Then  $\sig( A)=\sig_0(\ve)\cup\sig_1(\ve)$ is  a spectral decomposition of $\sig( A)$. Let $E= E_0(\ve)\oplus  E_1(\ve)$ be the corresponding direct sum  decomposition of  $ E$, where $E_1(\ve)=\oplus_{\mu\in\sig_1(\ve)}\mb{GE}_\mu(A)$ is  finite-dimensional (recall that $\mb{GE}_\mu(A)$ is  finite-dimensional for every $\mu\in\sig(A)\sm \sig_e(A)$).  A similar argument as in the verification of \eqref{e:3.0a} with some corresponding modifications applies to show that $ P\cap  E_1(\ve)\ne\{0\} .$

%\vs Pick a number  $\lam$ with $1/(r_e+\ve+\de)<\lam<1/(r_e+\ve)$. It is clear that
%\be\label{e:4.1}r_0:=\lam(r_e+\ve)<1,\hs\,r_1:=\lam(r_e+\ve+\de)>1.\ee
%Let $\~ A=\lam  A$. Then $\sig(\~ A)=\~\sig_0(\ve)\cup \~\sig_1(\ve)$, where $\~\sig_i(\ve)=\lam\sig_i(\ve)\,\,(i=0,1)$, and the corresponding spectral  decomposition of $ E$ remains the same as above. %, i.e., $E= E_0(\ve)\oplus  E_1(\ve)$.
%Noticing that  $$\~\sig_0(\ve)\subset\{z\in\bbC:\,\,\,|z|\leq r_0\},\hs \~\sig_1(\ve)\subset\{z\in\bbC:\,\,\,|z|\geq r_1\},$$ by \eqref{e:4.1} we deduce that  $$\mb{$||\~ A^ku||\ra 0$ ($u\in  E_0(\ve)$),\hs and $||\~ A^ku||\ra \8$  ($u\in  E_1(\ve)$)}$$ as $k\ra\8$.
%Now using some similar argument as in the verification of \eqref{e:3.0a}, one can easily check  that $  P\cap  E_1(\ve)\ne\{0\} $.

Let $A_1(\ve)=A|_{E_1(\ve)}$, and $P_1(\ve)=P\cap E_1(\ve)$. Then the operator $A_1(\ve)$ leaves $P_1(\ve)$ invariant.
  Choose an $\ve'\in(0,r-r_e)$ such that $r_e+\ve'$ is close to $r$ so that $\sig_1(\ve')\subset \bbS_r$; see Fig. \ref{fig:4.1} (b). Applying Theorem \ref{t:4.3} to $A_1(\ve')$ and  $ P_1(\ve')$,  one immediately concludes    that $r\in\sig(A_1(\ve'))=\sig_1(\ve')$ with a corresponding eigenvector $w\in P_1(\ve')$. Hence assertion (1) holds true.
\vs
In what follows we check the validity of assertions (2)-(4).\,

Let $\mu $ be an eigenvalue of $ A$ with  $|\mu|>r_e$. Pick a number  $\ve=\ve(\mu)>0$ sufficiently small so that $\mu\in \sig_1(\ve)$. Then clearly $r\in \sig_1(\ve)$.

For notational simplicity, we rewrite $E_1(\ve)=E_1$, $P_1(\ve)=P_1$, and $A_1(\ve)=A_1$.
 Suppose  that $  P$ contains a generalized eigenvector $\xi$ of $A$ pertaining to $\mu$. Since $\xi \in\mb{GE}_\mu(A)\subset E_1 $, we naturally have $\xi\in \mb{GE}_\mu(A_1 )$. Note that $\xi\in P_1 $. Because $r\in\sig(A_1)$ with a principal eigenvector $w\in P_1 $, Lemma \ref{l:4.5} applies to $ A_1 $. This allows us to deduce   that $\mu\in[0,r]$ with  a corresponding eigenvector $v\in  P_1 \subset P$, which completes the proof of assertion (2).

Now assume that $\mu\ne r$. Then by  Lemma \ref{l:4.5} we deduce that the eigenvectors of $ A_1 $ corresponding to  $\mu$ are contained in $ E_1 \sm\mb{int}_{E_1 }P_1 $. On the other hand, if $\o  P$ contains  an eigenvector $v $ of $\mu$, then  one easily verifies that $\mb{int}_{E_1 }P_1 \ne\emp$ and  $v \in \mb{int}_{E_1 }P_1 $. This leads to a contradiction and proves assertion (3).

             If $\o P\ne \emp$ and  contains a principal  eigenvector $u$, then  as above we have  $u\in \mb{int}_{E_1 } P_1 $. Assertion (4) then directly  follows from  Lemma \ref{l:4.5} (3). $\bx$

\bt\label{t:6.3}  Let  $A$ be  a positive operator with $r_e<r$. If  $ P$ is a solid cone and  $ A$  is rotationally strongly positive, then we have
  \benu
  \item[$(1)$] $r$ is an eigenvalue of $A$ with a  principal eigenvector $w\in P$;
\item[$(2)$] the algebraic and the geometric multiplicities  of  $r$ coincide;
 \item[$(3)$] $(\bbS_1 \xi)\,\cap \o P\ne\emp$ for every principal eigenvector $\xi$;
\item[$(4)$]  $\mb{\em GE}_\mu(A)\cap  P=\{0\}$ for any $\mu\in\sig(A)\sm\{r\}$ with $|\mu|>r_e$; and
\item[$(5)$] $|\mu|<r$ for all $\mu\in \sig(A)\sm\{ r\}$.
\eenu
\et
{\bf Proof.} Assertion (1) follows from Theorem \ref{t:6.1} (1). So we only need to verify assertions (2)-(5). We use the same notation as in the proof of Theorem \ref{t:6.1}.

Let $\mu $ be an eigenvalue of $ A$ with  $|\mu|>r_e$. Pick a number  $\ve=\ve(\mu)>0$ sufficiently small so that $\mu\in \sig_1(\ve)$. %Clearly $r\in \sig_1(\ve)$.
Let $E_1(\ve)=E_1$, and $A_1(\ve)=A_1$. Recall that $E_1$ is a finite-dimensional space and $P_1:=P\cap E_1\ne\{0\}$. Now we show that $\mb{int}_{E_1}P_1\ne\emp$ and $A_1$ is rotationally strongly positive
(with respect to $P_1$), and  the conclusions then  immediately follow from Theorem \ref{t:5.2}.

First, the same argument leading to \eqref{e:5.2} applies to show that there is a principal eigenvalue $\xi_1\in \o P$. As $\xi_1\in E_1$, one easily checks  that $\xi_1\in \mb{int}_{E_1}P_1$, that is, $\mb{int}_{E_1}P_1\ne\emp$.
Now let $u\in P_1\sm\{0\}$. Then by the rotational strong positivity of $A$, $\bbS_1(Au)\cap \o P\ne\emp$. As $Au\in E_1$ (and hence $\bbS_1(Au)\subset E_1$),  we have
$$\bbS_1(A_1u)\cap \mb{int}_{E_1}P_1\supset \bbS_1(Au)\cap \o P\ne\emp.$$ Hence $A_1$ is rotationally strongly positive.
$\bx$

\br{\em
In  \cite{Rugh} Rugh
generalized the notions of a real Birkhoff cone and its Hilbert's metric to complex spaces and obtained  a set  of spectral gap theorems. In particular, the classical PF theorem and  KR theorem were extended in some ways to complex matrices and operators.
Our approach in extending the PF theory to complex operators is of a pure dynamical nature, and is significantly different from the one in \cite{Rugh}.
}\er

 \section{The KR Type Theorems for Bounded Real Operators}\label{s:7}
% As direct consequences  of the results given in previous sections,  in this section we  present  some generalized versions of the classical KR theorem for real operators.

Now we focus  our attention on the real case.  So let $X$ be a   real Banach space, and  $P$ a cone in $X$.
 Let $A\in\sL(X)$.

Denote by
$\bbX=X+\mathbf{i} X$   the complexification of $X$ equipped with the norm
$$
\|u\|=\sup_{\theta\in[0,2\pi]}\|(\cos \theta) x+(\sin \theta) y\|,\Hs \A\,u=x+\bi y\in\bbX.
 $$
Then the space $X$ can be treated as a subspace of $\bbX$ because of the isometric embedding $j:X\ni x\mapsto x+\bi 0\in \bbX.$

There is a natural extension of the operator $A$ on $\bbX$, denoted by $\bbA$,
$$ \bbA u=Ax+\mathbf{i}Ay,\Hs \A\,u=x+\mathbf{i}y\in  \bbX.$$
$\bbA$ is called a complexification of $A$. Given   $\mu\in\sig(A):=\sig(\bbA)$, set
\be\label{e:GE}\ba{ll}
         \mb{GE}_\mu(A)=\{\mb{Re}\,\xi:\,\,\xi\in \mb{GE}_\mu(\bbA)\}=\{\mb{Im}\,\xi:\,\,\xi\in \mb{GE}_\mu(\bbA)\}\,.
         %\{\mb{\em Re}\,\xi:\,\,\xi\in \mb{GE}_\mu(\bbA)\}\Cup \{\mb{\em Im}\,\xi:\,\,\xi\in \mb{GE}_\mu(\bbA)\},
         \ea
         \ee
  (The second equality in \eqref{e:GE} is due to the fact that if $\xi\in \mb{GE}_\mu(\bbA)$, then $\pm \mathbf{i}\xi\in \mb{GE}_\mu(\bbA)$.)

    Note that by Remark \ref{r:2.9}, the cone $P$   is full if and only if it is total.

Let $\bbP= P+\mathbf{i} P$. Clearly  $\bbP$ is a cone in $\bbX$.
   We first give two basic lemmas concerning $\bbP$ that will be used in this section.
\bl\label{l:7.1}
If $P$ is a total cone   in $X$,  then $\bbP$ is full in $ \bbX$.
\el
{\bf Proof.} Let $ \bbX_0$ be a closed  subspace of $ \bbX$. Assume that $\bbP\subset \bbX_0$. Then we have  $P\subset \bbP\subset  \bbX_0$.
Set $X_0=\bbX_0\cap X$. Clearly $X_0$ is a closed real subspace of $X$ with $P\subset X_0$. Since $P$ is total  in $X$, we deduce  that $X_0=X$. Hence we see that $X\subset \bbX_0$. But since  $ \bbX_0$ is a complex Banach space, it then follows that  $\bbX=X+\mathbf{i}X\subset  \bbX_0$. Therefore $\bbX_0=\bbX$. $\bx$

\bl\label{l:rsp}  Let $P$ be a solid cone. If $A$ is strongly positive, then $\bbA$ is rotationally strongly positive (with respect to $\bbP$).
\el
{\bf Proof.} Let $u=x+\mathbf{i} y\in \bbP\sm\{0\}$. If $x\ne0\ne y$, then $Ax,Ay\in\o P$. Consequently $ \bbA u\in\o{\bbP}$. Now assume that $x=0\ne y$. Then by the  strong positivity of $A$ we have    $Ay\in\o P$. Taking  $z=(1-\mathbf{i})/\sqrt 2$, one finds   that
 $$
 z \bbA u=\frac{\sqrt2}{2}(1-\mathbf{i})(\mathbf{i} A y)=\frac{\sqrt2}{2}\(A y+\mathbf{i} Ay\)\in\o{\bbP}.
 $$

In the  case where  $x\ne 0=y$, in a similar manner as above  it can be shown that there is $z\in \bbS_1$ such that $z \bbA u\in \o{\bbP}$. $\bx$

%Now we state and prove the generalized versions of the KR theorem for real operators.
%Throughout this section we  always assume that %the operator  $A$ satisfies the following standing assumption: \vs
% (H) \,$r_e:=r_e(A)<r(A):=r$.

%%%%%%%%%%%%%%%%%%%%%%%%

\subsection{A general version of the KR theorem}
\vs

Let us first give  a general version of the  KR  Theorem.  For simplicity, we write $r_e(A)=r_e$, and $r(A)=r$.
\bt\label{t:7.4}   Let $A$ be a  positive operator with  $r_e<r$.  If  $P$ is total, then the following assertions hold:
        \benu
        \item[$(1)$] $r$ is an eigenvalue of $A$ with  a principal  eigenvector  $u\in P$.
\item[$(2)$] If $\o P\ne\emp$  and contains a principal eigenvector $v$ of $A$,  the algebraic and the geometric multiplicities  of $r$ coincide.
              \item[$(3)$] Let $\mu$ be an eigenvalue of $A$ with $|\mu|>r_e$. If $\mu\in\bbC\sm [0,\8)$, then $$\mb{\em GE}_\mu(A)\cap P=\{0\}.$$
                  \item[$(4)$] All eigenvectors of $A$ pertaining to other eigenvalues  $\mu\ne r$ with $|\mu|>r_e$ are contained in $X\sm\o P$.
\eenu
\et\vs
\br\label{r:7.1}{\em Assertion $(1)$ is already known in the literature for non-compact  operators; see e.g. \cite[Section 3, Theorem 1]{Edm} for  the case where  $P$ is reproducing (meaning that $X=P-P$),  and \cite[Theorem 1]{ZS} and \cite[Corollary 2.2]{Nuss} for the more general case where $P$ is total.

As far as we know, assertions (2)-(4) are new under the hypotheses in the theorem.}
\er

\noindent
{\bf Proof of Theorem \ref{t:7.4}.} (1) \,Lemma \ref{l:7.1} asserts that $\bbP$ is full in $\bbX$. Thus by Theorem \ref{t:6.1} we have  $r\in\sig(A)$; furthermore,    $ \bbA$  has  a corresponding eigenvector  $w\in \bbP$.
Let $w=u+\mathbf{i}v$, $u,v\in X$. We may assume $u\ne 0$. Then by the definition of $ \bbA$  it can be easily seen  that $u$ is an eigenvector of $A$. Note that $w\in \bbP$ implies $u\in P$.
       \vs
       (2) If $\o P$ contains a principal eigenvector $v$ of $A$, then $\xi=v+\mathbf{i}v\in\o\bbP$ and is a principal  eigenvector  of $\bbA$. The conclusion then follows from Theorem \ref{t:6.1}  (4).
       \vs
(3) \,%It is known that $\mb{GE}_\mu(A)$ is a finite dimensional space.
Suppose on the contrary that $P_0:=\mb{GE}_\mu(A)\cap P\ne\{0\}$. Then $P_0$ is a cone in the subspace $X_0=\mb{GE}_\mu(A)$. Let $A_0=A|_{X_0}$. Since $X_0$ is $A$-invariant, we trivially have   $A_0P_0\subset P_0$. We may assume that  $P_0$ is total in $X_0$; otherwise one can use the space $Y_0=\ol{P_0-P_0}$ in place of $X_0$. (Note that $Y_0$ is $A_0$-invariant.) By virtue of assertion (1) we deduce that $r(A_0)$ is an eigenvalue. This leads to a contradiction because $\sig(A_0)=\{\mu\}$.
\vs
           (4) \,If $ \o{P}$ contains an eigenvector $\eta$ of $A$ corresponding to a real  eigenvalue  $\mu\ne r$ with  $|\mu|>r_e$, then $\xi=\eta+\mathbf{i}\eta$ is an eigenvector  of $ \bbA$ corresponding  to $\mu$. Clearly  $\xi\in\o{\bbP}$. This contradicts Theorem \ref{t:6.1} (3) and verifies   assertion (4).
    $\bx$

\Vs
We infer from Theorem \ref{t:7.4} (4) that non-principal eigenvectors cannot occupy the interior $\o P$ of the cone. However, the following easy example indicates that  under the hypotheses of the theorem, they may lie on the boundary   $\pa P$.
  \Vs
\noindent{\bf Example 7.1.} Let $X=\R^2$, and $P=\R^2_+$. Then  the  matrix  $A=\(\begin{matrix}1&1\\ 0&2 \end{matrix}\)$ has two eigenvectors $w_1=(1,1)^T$ and $w_2=(1,0)^T$ (where $(x,y)^T$ denotes the transpose of $(x,y)$) corresponding to eigenvalues $r=2$ and $\mu=1$, respectively. Clearly $w_1\in\o P$, whereas $w_2\in \pa P$.

\subsection{Strong versions of the KR theorem}\vs
It is well known that the notions of strong positivity, irreducibility and primitivity  play crucial roles in the PF theory for real operators. A common feature of these notions  is that they exclude  the possibility that  eigenvectors of a  given  operator  may lie on the boundary of the cone.
 \bd Let $P$ be a solid cone. We say that $A$ is  weakly irreducible, if the boundary  $\pa P$ of $P$  contains no   eigenvectors of $A$ pertaining to  nonnegative  eigenvalues.
\ed

 It is almost obvious  that strongly positive operators and primitive operators are weakly irreducible.
 (Recall that  $A\in\sL(X)$ is  {\em primitive}, if there is an integer $m\geq 1$ such that $A^m(P\sm\{0\})\subset \o P$; see \cite[pp. 285]{LN}.)
In what follows we  show  that  irreducibility implies weak irreducibility    in the above terminology. Thus we assume that $A\in\sL(X)$ is {\em irreducible}, i.e., there exists  $\lam>r( A)$ such that
\be\label{e:5.1}R_\lam(A)( P\sm\{0\})\subset \o P.\ee
   We argue by contradiction and suppose that $A$ were   not weakly irreducible. Then $A$ would have  an eigenvector   $v$ pertaining to a real eigenvalue $\mu\geq 0$ such that  $v\in\pa P$. Observing  that $(\lam- A)v=(\lam-\mu) v$, where $\lam$ is  the number in \eqref{e:5.1},  we deduce that   $(\lam -\mu)R_\lam(A)v=v\in\pa P$. Since $\lam-\mu>0$, it follows  that $R_\lam(A)v\in\pa P$, which contradicts  \eqref{e:5.1}.

\br{\em  In the finite-dimensional case an equivalent  definition of  irreducibility  for linear operators can be found in \cite{LN} by using the notion  of  {\em faces} of  cones; see  \cite[Definition 1.1.4, Proposition 1.1.5]{LN} for details.  We also refer the interested  reader to \cite[pp. 341]{Schab} or \cite[pp. 149]{de} for the definition of irreducibility  on Banach lattices which makes use of the concept of {\em ideals}.
}
\er
\vs
As a simple  consequence of Theorem \ref{t:7.4}, one immediately obtains a refined version of the KR theorem for weakly irreducible operators.

\bt\label{t:7.5} Let $P$ be a solid cone, and   $A$ a positive operator with  $r_e<r$. Suppose that  $A$ is weakly irreducible. Then
                \benu
        \item[$(1)$] $r$ is a  simple  eigenvalue of $A$ with  a principal  eigenvector  $w\in \o P$; and
                    \item[$(2)$]
                     $\mb{\em GE}_\mu(A)\cap P=\{0\}$ for any  $\mu\in\sig(A)\sm \{r\}$ with $|\mu|>r_e$.
           \eenu
\et
{\bf Proof.} (1) By virtue of assertion (1) in Theorem \ref{t:7.4} we conclude that  $r$ is an eigenvalue of $A$ with  a corresponding eigenvector  $w\in P$.  The weak irreducibility of $A$ then implies that $w\in\o P$. Thus  by  Theorem \ref{t:7.4} (2), we deduce that  $r$ has the same  algebraic and  geometric multiplicities.

Now suppose  that $A$ has another principle  eigenvector $v$  which is not a multiple of $w$. Let $\pi=\{s w+tv:\,\,s,t\in\R\}$ be the real plane spanned by $w$ and $v$.  Since $w\in\o P$, it is easy to see that $\pi\cap \pa P$ contains a half-line. Take a $u\in \pi\cap \pa P $, $u\ne0$. Then  $u$ is a principal eigenvector of $A$, contradicting the weak  irreducibility of $A$. This finishes the proof of assertion (1).
\vs
(2) \,In view  of Theorem \ref{t:7.4} (3), one only needs to consider the case where $\mu$ is a real eigenvalue of $A$.

We argue by contradiction and suppose that $\mb{GE}_\mu(A)\cap P\ne\{0\}$. Then as in the proof of Theorem \ref{t:7.4} (3), one deduces that $P_0=\mb{GE}_\mu(A)\cap P$ is a cone in $X_0=\mb{GE}_\mu(A)$ which is left invariant by $A_0:=A|_{X_0}$. It can be  assumed that  $P_0$ is total in $X_0$.  Thus by Theorem \ref{t:7.4} (1), $A_0$ has an eigenvector $u\in P_0$ pertaining to $\mu$. Further by Theorem \ref{t:7.4} (4),  we conclude that $u\in \pa P$, which contradicts the weak irreducibility of $A$.  $\bx$

\br{\em For a non-compact  operator $A$ as considered in Theorem \ref{t:7.5},  in  \cite{ZhangL} Zhang  proved that the principal eigenvalue $r$ is simple under the hypothesis of strong positivity on $A$; see \cite[Theorem 1.3]{ZhangL}. Here we have proved  the same conclusion  under a more general  and weaker hypothesis.

To the best of our knowledge,  Theorem \ref{t:7.5}  (2)   is new  even if in the case where $A$ is assumed to be strongly positive as considered in \cite{ZhangL}.
}\er

The following easy example indicates that in general, the peripheral spectrum of  a weakly irreducible operator   $A$ as considered in Theorem \ref{t:7.5} may contain  more than one eigenvalues.
\Vs
\noindent{\bf Example 7.2.} Let $X=\R^2$,  $P=\R^2_+$, and   $A=\(\begin{matrix}0&1\\[-1ex] 1&0 \end{matrix}\)$.
Then  one easily  verifies that  $R_\lam(A)=\frac{1}{3}\(\begin{matrix}2&1\\[-1ex] 1&2 \end{matrix}\)$ for $\lam=2$, and hence  $R_\lam(A)( P\sm\{0\})\subset \o P$. Therefore  $A$ is weakly irreducible. However, we have $\sig(A)=\{\pm 1\}$.
\Vs

 Theorem \ref{t:7.8} below  is a modified version of \cite[Theorem 1.3]{ZhangL}, which asserts that if $A$ is strongly positive then its peripheral spectrum consists of exactly the principal eigenvalue $r$. Besides we have added in a new  conclusion (i.e. assertion (2)) to \cite[Theorem 1.3]{ZhangL}, the real point here is that we have given  an essentially   different dynamical proof for the theorem.

\bt\label{t:7.8} Let  $P$ be a solid  cone. Suppose that  $A$ is strongly positive and that   $r_e<r$. Then
 \benu
        \item[$(1)$] $r$ is a  simple  eigenvalue of $A$ with  a principal  eigenvector  $w\in \o P$;
                    \item[$(2)$]
                     $\mb{\em GE}_\mu(A)\cap P=\{0\}$ for any  $\mu\in\sig(A)\sm \{r\}$ with $|\mu|>r_e$; and
                     \item[$(3)$] $|\mu|<r$ for all $\mu\in\sig(A)\sm \{r\}$.
           \eenu
\et

\noindent{\bf Proof.} Since the strong positivity of $A$ implies that $A$ is weakly irreducible, the validity of assertions (1) and (2) directly follows from Theorem \ref{t:7.5}. Assertion (3) is a simple   consequence of Theorem \ref{t:6.3} and Lemma \ref{l:rsp}. $\bx$

\section*{Acknowledgement}

We are indebted to an  anonymous referee who has  checked the manuscript very carefully for several times. He/She has pointed out several  gaps in  the proofs of some fundamental results in the original version of this manuscript and given us  many constructive  comments and suggestions, which  greatly helped us to improve  the quality of the manuscript.

Our sincere thanks also go to Professor  R. Nagel  for  bringing to our attention some important works in the literature and helpful discussions on   the Perron-Frobenius theory.

\section*{References}
{\small

\end{document}